\documentclass[1p]{elsarticle}
\usepackage{amsmath}
\usepackage{amssymb}
\usepackage{amsthm}
\usepackage{lineno}
\usepackage{subfigure}
\usepackage{graphicx}
\usepackage{color}
\usepackage{xspace}
\usepackage{hyperref} 
\graphicspath{{../MAfigures/}}

\newcommand{\mao}{MA}
\newcommand{\delu}{v}
\newcommand{\trace}{\text{trace}}

\newcommand{\bq}{\begin{equation}}
\newcommand{\eq}{\end{equation}}
\newcommand{\R}{\mathbb{R}}
\newcommand{\dm}{d}
\newcommand{\Rd}{\R^\dm}
\newcommand{\Dt}{\mathcal{D}}
\newcommand{\abs}[1]{\left\vert#1\right\vert}
\newcommand{\norm}[1]{\left\vert#1\right\vert}
\newcommand{\DD}[2]{\frac{\partial^2 #1}{\partial {#2}^2}}
\newcommand{\MA}{{Monge-Amp\`ere}\xspace}

\newcommand{\blue}[1]{{#1}}

\newcommand{\grad}{\nabla}
\newcommand{\e}{\epsilon}
\newcommand{\ex}[1]{\times 10^{- #1}}
\newcommand{\G}{\mathcal{G}}
\newcommand{\bO}{\mathcal{O}}

\newcommand{\xv}{\mathbf{x}}
\newcommand{\xo}{\mathbf{x}_0}

\newtheorem{theorem}{Theorem}
\newtheorem{lemma}[theorem]{Lemma}

\newdefinition{remark}{Remark}
\newdefinition{definition}{Definition}
\newdefinition{example}{Example}

\begin{document}

\begin{frontmatter}

\title{Fast finite difference solvers for singular solutions of the elliptic Monge-Amp\`ere equation}

\author[sfu]{B. D. Froese}
\ead{bdf1@sfu.ca}
\ead[url]{http://www.divbyzero.ca/froese/}
\author[sfu]{A. M. Oberman\corref{cor1}}
\ead{aoberman@sfu.ca}
\ead[url]{http://math.sfu.ca/~aoberman}
\cortext[cor1]{Corresponding Author}
\address[sfu]{Department of Mathematics, Simon Fraser University\\
Burnaby, British Columbia, Canada, V5A 1S6}

\begin{abstract}
The elliptic Monge-Amp\`ere equation is a fully nonlinear Partial Differential Equation 
which originated in geometric surface theory, and has been applied in dynamic meteorology, elasticity, geometric optics, image processing and image registration.  
Solutions can be singular, in which case standard numerical approaches fail.

In this article we build a finite difference solver for the \MA equation, which converges even for singular solutions.
Regularity results are used to select \emph{a priori} between a stable, provably convergent monotone discretization and an accurate finite difference discretization in different regions of the computational domain.  This allows singular solutions to be computed using a stable method, and regular solutions to be computed more accurately.
The resulting nonlinear equations are then solved by Newton's method.

Computational results in two and three dimensions validate the claims of accuracy and solution speed. 
 A computational example is presented which demonstrates the necessity of the use of the monotone scheme near singularities.
\end{abstract}

\begin{keyword}
Fully Nonlinear Elliptic Partial Differential Equations \sep
 Monge Amp\`ere equations \sep 
 Nonlinear Finite Difference Methods\sep 
 Viscosity Solutions\sep 
 Monotone Schemes\sep 
 Convexity Constraints
\end{keyword}

\end{frontmatter}
\linenumbers

\section{Introduction}\label{sec:intro}
In this article we build a finite difference solver for the \MA equation, which converges even for singular solutions.
Regularity results are used to select \emph{a priori} between two discretizations in different regions of the computational domain.  
Near possible singularities, a stable, provably convergent monotone discretization is used. 
Elsewhere a more accurate discretization is used.
This allows singular solutions to be computed using a stable method, and regular solutions to be computed more accurately.
The resulting nonlinear equations are then solved by Newton's method, which is fast, $\bO(M^{1.3})$, where $M$ is the number of data points, independent of the regularity of the solution.




\subsection{The setting for equation}

The \MA equation is a fully nonlinear Partial Differential Equation (PDE).  
\bq
\label{MA}\tag{MA}
\det(D^2u(x)) = f(x), \quad \text{for $x$ in }\Omega.
\eq
The \MA operator, $\det(D^2u)$,  is the determinant of the Hessian of the function $u$.   The equation is augmented by the convexity constraint
\bq\label{convex}\tag{C}
u \text{ is convex, }
\eq
which is necessary for the equation to be elliptic.  
The convexity constraint is made explicit for emphasis: it is necessary for uniqueness of solutions and it is essential for numerical stability.  

While other boundary conditions appear naturally in applications, we consider the simplest boundary conditions: the Dirichlet problem in a convex bounded subset~$ \Omega \subset \Rd$ 
with boundary~$\partial\Omega$,
\bq\tag{D}\label{Dirichlet}
u(x) = g(x), \quad \text{for $x$ on }\partial\Omega.
\eq
Under suitable assumptions on the domain and the functions $f(x), g(x)$, recalled in~\autoref{sec:regularity}, there exist unique classical ($C^2$) solutions to~\eqref{MA},~\eqref{convex}.  
However, when these conditions fail, solutions can be singular.  
For singular solutions, the correct notion of weak solutions must be used.
In this case, novel discretizations and solutions methods must be used to approximate the solution.

\subsection{Applications}
The PDE~\eqref{MA} is a geometric equation, which goes back to Monge and Amp\`ere (see~\cite{EvansSurvey}).   
The equation naturally arises in geometric problems of existence and uniqueness of surfaces with proscribed metrics or curvatures~\cite{Bakelman, PogBook}.  
Early applications identified in~\cite{olikerprussner88} include dynamic meteorology, elasticity, and geometric optics~\cite{Haltiner, Kasahara, Stoker, Westcott}.  
For an application of \MA equations in mathematical finance, see~\cite{StojanovicMAapp}.

The \MA equation arises as the optimality conditions for the problem of optimal mass transport with quadratic cost~\cite{EvansSurvey, Ambrosio,Villani}.
This application of the \MA equation has been used in many areas:
image registration~\cite{Haker, HakerRegistration, RehmanRegistration},
 mesh generation~\cite{Delzanno, DelzannoGrid, Budd}, 
 reflector design~\cite{GlimmOlikerReflectorDesign}, 
 and astrophysics (estimating the shape of the early universe)~\cite{FrischUniv}.

\newcommand{\fv}{\mathbf{g}}

The problem here is to find a mapping $\fv(x)$ that moves the measure $\mu_1(x)$ to $\mu_2(y)$ and minimizes the transportation cost functional
\[ 
\int_{\R^d} \norm{ x-\fv(x) }^2 \,d\mu_1. 
\]
The optimal mapping is given by $\fv = \grad u$,  where $u$ satisfies the \MA equation
\[ \det(D^2u(x)) = \mu_1(x)/\mu_2(\nabla u(x)). \]  
In this situation, the Dirichlet boundary condition~\eqref{Dirichlet} is replaced by the implicit boundary condition
\bq\label{BC2}
\fv(\cdot) : \Omega_1 \to \Omega_2
\eq
where the sets $\Omega_1$ and $\Omega_2$ are the support of the measures~$\mu_1,\mu_2$.
These boundary conditions are difficult to implement numerically; we are not aware of an implementation using PDE methods.  For many applications, 
 both domains are squares, and a simplifying assumption that edges are mapped to edges allows {Neumann} boundary conditions to be used.  In other applications, periodic boundary conditions are used.

In other problems, the \MA operator appears in an \emph{inequality constraint} in a variational problem for optimal mappings where the cost is no longer the transportation cost.  Here the operator has the effect of restricting the local area change on the set of admissible mappings, see~\cite{HaberTransport} or~\cite{CohenOr}.  

\subsection{Related numerical works}
Despite the number of applications, until recently there have been few numerical publications devoted to solving the \MA equation.  We make a distinction between numerical approaches with optimal transportation type boundary conditions~\eqref{BC2} and the standard Dirichlet boundary conditions~\eqref{Dirichlet}.  
In the latter case, a number of numerical methods have been recently proposed for the solution of the \MA equation.

An early work is~\cite{olikerprussner88}, which presents a discretization which converges to the Aleksandrov solution in two dimensions.  
Another early work by Benamou and Brenier~\cite{BenBren} used a fluid mechanical approach to compute the solution to the optimal transportation problem.

For the problem with Dirichlet boundary conditions which is treated here, a series of papers have recently appeared by two groups of authors, Dean and Glowinski~\cite{DGnum2008, DGaug, GlowinksiICIAM}, and Feng and Neilan,~\cite{FengMA, FengFully}.  The methods introduced by these authors perform best in the regular case and can break down in the singular case.  See~\cite{BeFrObMA} a more complete discussion.

We also mention the works~\cite{LoeperMA}, in the periodic case, and~\cite{Delzanno} for applications to mappings.  The method of~\cite{MAFrisch} treats the problem of periodic boundary conditions in odd dimensional space.

\subsection{Numerical challenges}

When the conditions for regularity are satisfied, classical solutions can be approximated successfully using a range of standard techniques (see, for example works such as~\cite{DGnum2008, DGaug, GlowinksiICIAM}, and~\cite{FengMA, FengFully}).  
However, for singular solutions, standard numerical methods break down: either by becoming unstable, poorly conditioned, or by selecting the wrong (non-convex) solution.

\subsubsection*{Weak solutions}
For singular solutions, the appropriate notion of weak (viscosity or Aleksandrov) solutions must be used.  
Numerical methods have been developed which capture weak solutions:
Oliker and Prussner, in~\cite{olikerprussner88},  presents a method  which converges to the Aleksandrov solution.  
One of the authors introduced a finite difference method which converges to the viscosity solution in~\cite{ObermanEigenvalues}. Both of these methods were restricted to two dimensions.    However, methods which are provably convergent may have lower accuracy or slower solution methods than other methods which are effective for regular solutions.
In~\cite{ABtheory} we introduced a monotone discretization which is valid in arbitrary dimensions.  A proof of convergence to viscosity solutions is provided, as well as a proof of convergence of Newton's method.   

\subsubsection*{Convexity}
The convexity constraint is necessary for both uniqueness and stability.
In particular, the equation~\eqref{MA} fails to be elliptic if $u$ is non-convex (see~\autoref{sec:ellipt}).   
so instabilities can arise if the convexity condition~\eqref{convex} is violated, as demonstrated in~\autoref{sec:newtonfails}.
Any approximation of~\eqref{MA} requires some selection principle to choose the convex solution.  
This selection principle can be built in to the discretization, as in~\cite{ObermanEigenvalues}, or built in to the solution method, as in~\cite{BeFrObMA}.

\subsubsection*{Accuracy}
The convergent monotone schemes of~\cite{ObermanEigenvalues} and~\cite{ABtheory} use a wide stencil, and the accuracy of the scheme depends on the \emph{directional resolution}, which depends on the width of the stencil.  As we demonstrate below, for highly singular solutions, such as~\eqref{eq:cone}, the directional resolution error can dominate.  On the other hand, more accurate discretizations, such as standard finite differences, can be unstable for singular solutions.

\subsubsection*{Fast solvers}
Previous work by the authors and a coauthor~\cite{BeFrObMA} investigated fast solvers for~\eqref{MA}.  An explicit method was presented which was moderately fast, independent of the solution time.  
For regular solutions, a faster (by an order of magnitude) semi-implicit solution method was introduced (see~\autoref{sec:SI})  but this method was slower (by an order of magnitude) on singular solutions.

\section{Analysis and weak solutions}\label{sec:weak}
In this section we present regularity results and background analysis which inform the numerical approach taken in this work.  In particular, the regularity results of~\autoref{sec:regularity} are used to determine the discretization used in~\autoref{sec:hybrid}.  

The definition of viscosity solutions  and Aleksandrov solutions presented in \autoref{sec:viscosity}-\ref{sec:alex} are used to make sense of the weak solutions~\eqref{eq:c1} and~\eqref{eq:cone}, respectively.

\subsection{Regularity}\label{sec:regularity}
Under the following conditions, the \MA equation is guaranteed to have a unique $C^{2,\alpha}$ solution
Regularity results for the \MA equation have been established in~\cite{CafNirSpruck,UrbasDirichletMA,CaffEstimates}. We refer to the book~\cite{Gutierrez} for the following result.

\bq\label{conditions}
\begin{cases}
\text{The domain $\Omega$ is strictly convex with boundary $\partial\Omega\in C^{2,\alpha}$.}\\
\text{The boundary values $g \in C^{2,\alpha}(\partial\Omega)$.}\\
\text{The function $f \in C^\alpha(\Omega)$ is strictly positive.}
\end{cases}
\eq

\begin{remark}  In the extreme case, with $f(x) = 0$ for all $x \in \Omega$, the equation~\eqref{MA},\eqref{convex} reduces to the computation of the convex envelope of the boundary conditions~\cite{ObermanCENumerics, ObSilvConvEnv}.  In this case, solutions may not even be continuous up to the boundary and can also 
be non-differentiable in the interior.
\end{remark}

\begin{remark}
While is it usual to perform numerical solutions on a rectagle, 
regularity can break down in particular convex polygons~\cite{Villani,Pogorelov}.  
\end{remark}

\subsection{Viscosity solutions}\label{sec:viscosity}
We recall the definition of viscosity solutions~\cite{CIL}, which are defined for the \MA equation in~\cite{Gutierrez}. 
\begin{definition}
Let $u \in C(\Omega)$ be convex and $f\geq0$ be continuous.  The function $u$ is a \emph{viscosity subsolution (supersolution)} of the \MA equation in $\Omega$ if whenever convex $\phi\in C^2(\Omega)$ and $x_0\in\Omega$ are such that $(u-\phi)(x)\leq(\geq)(u-\phi)(x_0)$ for all $x$ in a neighbourhood of $x_0$, then we must have
\[ \det(D^2\phi(x_0)) \geq(\leq)f(x_0). \]
The function $u$ is a \emph{viscosity solution} if it is both a viscosity subsolution and supersolution.
\end{definition}

\begin{example}[Viscosity solution of \MA]
We consider an  example which will later be solved numerically in two and three dimensions (sections~\ref{sec:2d}-\ref{sec:3d}).  
Consider~\eqref{MA} with solution and $f$ given by
\[ 
u(\xv) = \frac{1}{2}((\abs{\xv}-1)^+)^2, 
\qquad
f(\xv) = (1-1/\abs{\xv})^+. \]

This function, 
although it is not a classical $C^2$ solution of the \MA equation, is a viscosity solution.
\end{example}

\subsection{Aleksandrov solutions}\label{sec:alex}
Next we turn our attention to the Aleksandrov solution, which is a more general weak solution than the viscosity solutions.
Here $f$ is generally a measure~\cite{Gutierrez}.  We begin by recalling the definition of the normal mapping or subdifferential of a function.

\begin{definition}
The \emph{normal mapping} (\emph{subdifferential}) of a function $u$ is the set-valued function $\partial u$ defined by
\[ \partial u(x_0) = \{p  : u(x) \geq u(x_0) + p\cdot (x-x_0)\},\quad \text{ for all } x\in\Omega. \]
For a set $V \subset\Omega$, we define $\partial u(V) = \bigcup\limits_{x\in V}\partial u(x)$.
\end{definition}

Now we want to look at a measure generated by the \MA operator.  
\begin{definition}
Given a function $u\in C(\Omega)$, the \emph{\MA measure} associated with $u$ is defined as
\[ \mu(V) = \abs{\partial u(V)}\]
for any set $V \subset\Omega$.
\end{definition}

This measure naturally leads to the notion of the generalized or Aleksandrov solution of the \MA equation.
\begin{definition}
Let $\mu$ be a Borel measure defined in a convex set $\Omega\in\R^d$.  Then the convex function $u$ is an \emph{Aleksandrov solution} of the \MA equation
\[ \det(D^2u) = \mu \]
if the \MA measure associated with $u$ is equal to the given meaure $\mu$.
\end{definition}

\begin{example}[Aleksandrov solution]
As an example, we consider the cone and the the scaled Dirac measure
\[ u(\xv) = \norm{\xv}, 
\qquad
 \mu(V) = \pi \int_V \delta(\xv)\,d\xv. 
 \]
\end{example}

\subsection{A PDE for convexity}
The convexity constraint~\eqref{convex} is necessary for uniqueness, since without it, $-u$ is also a solution of~\eqref{MA}.

For a twice continuously differentiable function $u$, the convexity restriction~\eqref{convex} can be written as $D^2u$ is positive definite. 
Since we wish to work with less regular solutions,~\eqref{convex} can be enforced by the equation
\[
\lambda_1(D^2u) \ge 0,
\] 
understood in the viscosity sense~\cite{ObermanCENumerics, ObSilvConvEnv}, 
where $\lambda_1[D^2u]$ is the smallest eigenvalue of the Hessian of $u$.

The convexity constraint can be absorbed into the operator by
 defining
\bq\label{detplus}
{\det}^+(M) = \prod\limits_{j=1}^d\lambda_j^+ 
\eq
where $M$ is a symmetric matrix, with eigenvalues, $\lambda_1 \le \dots, \le \lambda_n$ and
\[
x^+ = \max(x, 0).
\] 
Using this notation,~\eqref{MA},\eqref{convex} becomes
\bq
\label{MAplus}\tag*{$(MA)^+$}
{\det}^+(D^2u(x))= f(x), \quad \text{for $x$ in }\Omega
\eq

\begin{remark}
Notice that there is a trade off in defining~\eqref{detplus}: the constraint~\eqref{convex} is eliminated but the operator becomes non-differentiable near singular matrices.
\end{remark}

\subsection{Linearization and ellipticity}\label{sec:ellipt}
The linearization of the determinant is given by
\[
\grad \det(M) \cdot N = 
\trace 
\left(   
	M_{adj} N	
\right )
\]
Where $M_{adj}$ is the adjugate~\cite{Strang}, which is the transpose of the cofactor matrix.  The adjugate matrix is positive definite if and only if $M$ is positive definite.
When the matrix $M$ is invertible, the adjugate, $M_{adj}$, satisfies
\bq\label{adjinv}
M_{adj} = \det(M) M^{-1}
\eq

We now apply these considerations to the linearization of the \MA operator.
When $u \in C^2$ we can linearize this operator as
\bq\label{eq:lin} 
\blue{
\grad_{M} \det(D^2u) \cdot v   = {\trace}\left( 
(D^2u)_{\blue{adj}}D^2v
\right ). 
}
\eq

\begin{example}
In two dimensions we obtain
\[ 
\grad_M \det(D^2u) v = 
u_{xx} v_{yy} + u_{yy} v_{xx} - 2 u_{xy} v_{xy}
\]
which is homogenous of order one in $D^2u$.
In dimension $d \ge 2$, the linearization is homogeneous order $d-1$ in $D^2u$.
\end{example}

The linear operator
\[
L[u] \equiv \trace{ A(x) D^2 u }
\]
is \emph{elliptic} if the coefficient matrix $A(x)$ is positive definite.

\begin{lemma}\label{lem:linell}
Let $u \in C^2$.  The linearization of the \MA operator,~\eqref{eq:lin} is elliptic if $D^2u$ is positive definite or, equivalently, if $u$ is (strictly) convex.
\end{lemma}

\begin{remark}
When the function $u$ fails to be strictly convex, the linearization can be degenerate elliptic, which affects the conditioning of the linear system~\eqref{eq:lin}.  When the function $u$ is nonconvex, the linear system can be ustable.
\end{remark}

The definition of a nonlinear elliptic PDE operator generalizes the definition of linear elliptic operator.  It also allows for the operators to be non-differentiable. 
\begin{definition}
Let the PDE operator $F(M)$ be a continuous  function defined on symmetric matrices.
Then $F(M)$ is \emph{elliptic} if it satisfies the monotonicity condition
\[ 
F(M) \leq F(N) \text{ whenever } M\leq N, 
\]
where for symmetric matrices $M\le N$ means $x^T M x \le x^T N x$ for all $x$.  
\end{definition}

\begin{example}
The operator ${\det}^+(M)$ is a non-decreasing function of the eigenvalues, so it is elliptic.
\end{example}

\section{The standard finite difference discretization}
We begin by considering the standard finite difference discretization of the \MA equation.  For brevity, we describe the discretization in two dimensions, but this is easily generalized to higher dimensions. 

This discretization does not enforce the convexity condition~\eqref{convex}, so it can lead to instabilities.  In particular, we show in~\autoref{sec:newtonfails} that Newton's method can become unstable if this discretization is used.

The \MA operator has a particularly simple form in two dimensions:
\[
\det(D^2u) = \DD{u}{x} \DD{u}{y}     - \left(\frac{\partial^2 u}{\partial x \partial y}\right)^2, \quad
\text{ in } \Omega \subset \R^2.
\]
In two dimensions, the natural discretization of the operator is given by
\bq\label{MANat}\tag*{$(MA)^{N}$}
\mao^{N}[u] \equiv (\Dt_{xx}u)(\Dt_{yy}u) - (\Dt_{xy}u)^2 
\eq
where, writing $h$ for the spatial resolution of the grid,
\begin{align*}
[\Dt_{xx}u]_{ij} &= \frac1 {h^2} 
\left(
{u_{i+1,j}+u_{i-1,j}-2u_{i,j}}
\right)
\\
[\Dt_{yy}u]_{ij} &= \frac{1}{h^2}
\left(
u_{i,j+1}+u_{i,j-1}-2u_{i,j}
\right)
\\
[\Dt_{xy}u]_{ij} &= \frac{1}{4h^2}
\left(
u_{i+1,j+1}+u_{i-1,j-1}-u_{i-1,j+1}-u_{i+1,j-1}
\right).
\end{align*}

\begin{remark}
There is no reason to assume that the standard discretization converges.
In fact, the two dimensional scheme has multiple solutions.  In~\cite{BeFrObMA} this discretization was used, but the the solvers were designed to select the convex solution. 
\end{remark}

\section{Convergent monotone discretization}\label{sec:monotone}
The method of~\cite{ObermanEigenvalues} describes a discretization of the two-dimensional \MA equation that converges to the viscosity solution. 
In~\cite{ABtheory} we introduced another discretization, which generalized to higher dimensions, and also converged to the viscosity solution.
Both methods require the use of a wide stencil scheme, which has an additional discretization parameter, the \emph{directional resolution}, explained below.

In addition to being monotone, which means it is provably convergent, the latter method discretizes the convexified version of the equation,~\ref{MAplus}, which is enough to ensure convergence of Newton's method.  The proof of this result can be found in~\cite{ABtheory}.

In this section we  present the convergent discretization, which will be used to build the hybrid solver.

\subsection{Wide stencils}
When we discretize the operator on a finite difference grid, we approximate the second derivatives by centred finite differences (spatial discretization).
In addition, we consider a finite number of possible directions $\nu$ that lie on the grid (directional discretization).

\begin{figure}
  \centering
  \subfigure[In the interior.]{
  \includegraphics[height=0.4\textwidth]
  {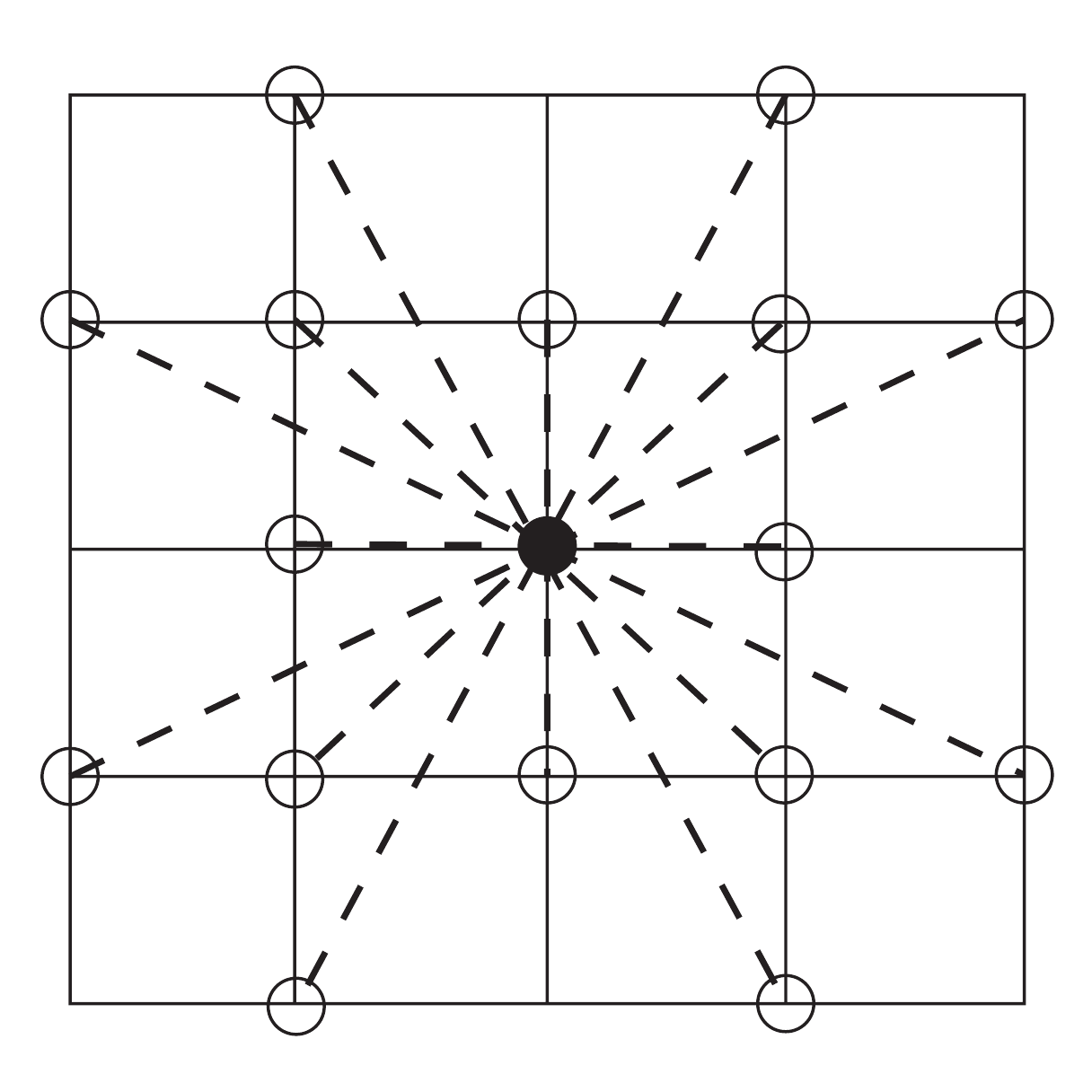}
  }                
  \subfigure[Near the boundary.]{
  \includegraphics[height=0.4\textwidth]
  {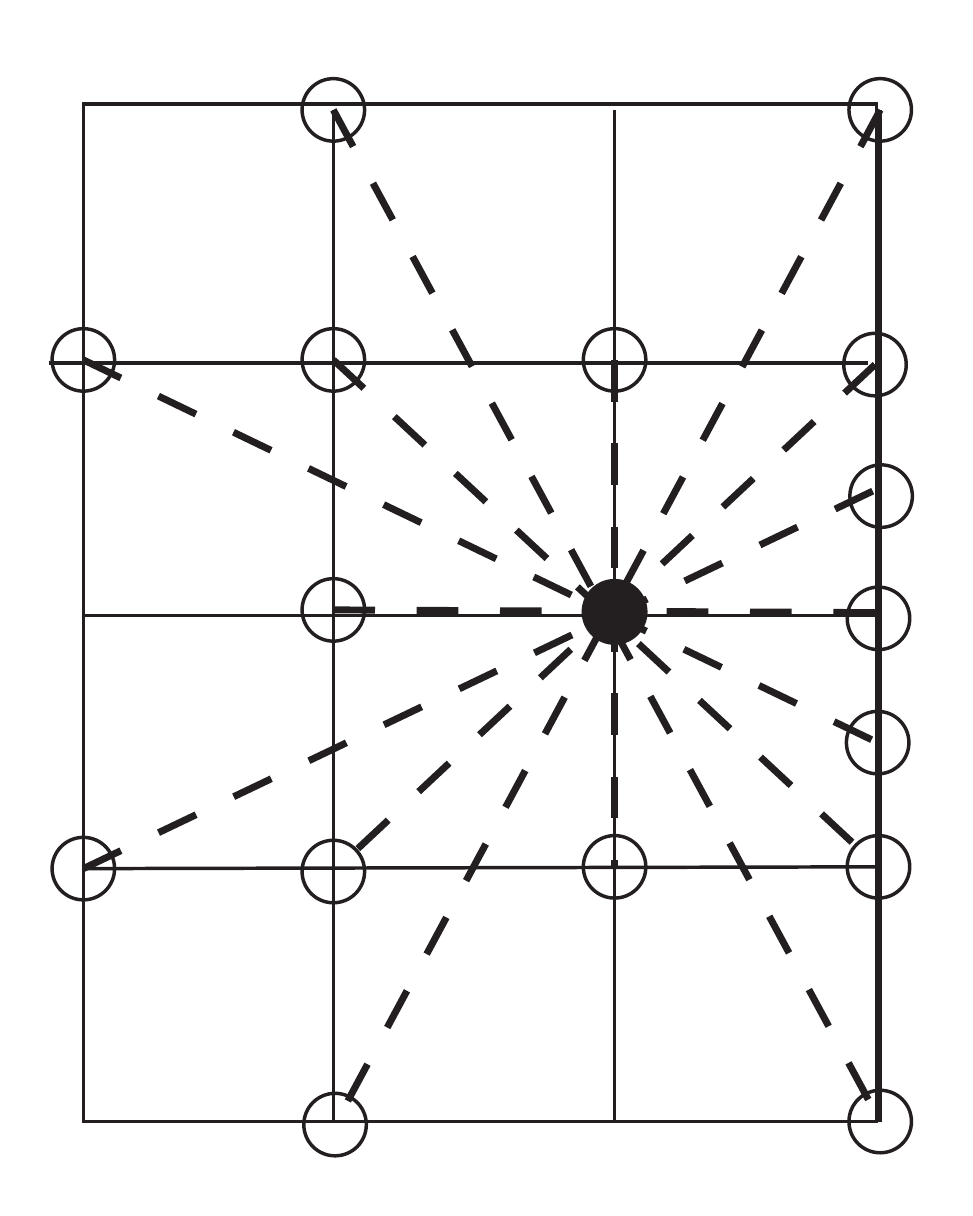}
  }
\caption{
Wide stencils on a two dimensional grid.
}
\label{fig:stencil}
\end{figure}

We consider the finite difference operator for the second directional derivative in the direction $\nu$, which lies on the finite difference grid.  
These directional derivatives are discretized by simply using finite differences on the grid
\[ \Dt_{\nu\nu}u_i = \frac{1}{\abs{\nu}h^2}\left(u(x_i + \nu h) + u(x_i - \nu h) - 2u(x_i)\right). \]
Depending on the direction of the vector $\nu$, this may involve a wide stencil.  At points near the boundary of the domain, some values required by the wide stencil will not be available; see~\autoref{fig:stencil}.  In these cases, we use interpolation at the boundary to construct a (lower accuracy) stencil for the second directional derivative; see~\cite{ObermanEigenvalues} for more details.

Since the discretization considers only a finite number of directions $\nu$, there will be an additional term in the consistency error coming from the angular resolution $d\theta$ of the stencil.  This angular resolution will decrease and approach zero as the stencil width is increased.

\subsection{Discretization of the convexified \MA operator}
In two dimensions, the largest and smallest eigenvalues of a symmetric matrix can be represented by the variational formula
\[
\lambda_1[A] = \min_{\norm{\nu} = 1} \nu^T A \nu, 
\qquad 
\lambda_2[A] = \max_{\norm{\nu} = 1} \nu^T A \nu.
\]
This formula was used in~\cite{ObermanEigenvalues} to build a monotone scheme for the~\MA operator, which is the product of the eigenvalues of the Hessian, by replacing the $\min, \max$ over all directions, by a finite number of grid directions.

In higher dimensions, the formula above does not generalize naturally.  Instead, in~\cite{ABtheory}, 
we used another characterization, which applied to positive definite matrices.

\begin{lemma}[Variational characterization of the determinant]\label{thm:eigs}
Let $A$ be a $d \times d$ symmetric positive definite matrix with eigenvalues $\lambda_j$ and let $V$ be the set of all  orthonormal bases of $\R^d$:
\[ V = \{(\nu_1,\ldots,\nu_d)  \mid \nu_j\in\R^d,\nu_i\perp\nu_j \text{ if }i\neq j, \|\nu_j\|_2 = 1\}. \]
Then the determinant of $A$ is equivalent to
\[ \prod\limits_{j=1}^d \lambda_j = \min\limits_{(\nu_1,\ldots,\nu_d)\in V} \prod\limits_{j=1}^d \nu_j^T A \nu_j. \]
\end{lemma}

We use Lemma~\ref{thm:eigs} to characterize the determinant of the Hessian of a convex $C^2$ function $\phi$ in terms of second directional derivatives of $\phi$.
\[ 
\det(D^2\phi) =  \min\limits_{(\nu_1,\ldots,\nu_d)\in V} \prod\limits_{j=1}^d \nu_j^T D^2\phi \nu_j
  = \min\limits_{(\nu_1,\ldots,\nu_d)\in V} \prod\limits_{j=1}^d \DD{\phi}{\nu_j}.
\]
The convexified \MA operator~\ref{MAplus} can then be represented by simply enforcing positivity of the eigenvalues, which leads to the following,
\[ 
{\det}^+(D^2 \phi) =  \min\limits_{\{\nu_1\ldots\nu_d\}\in V}
\prod\limits_{j=1}^{d} 
\left (\frac{\partial^2 \phi}{\partial\nu_j^2}\right )^+
.
\]
To discretize the operator on a finite difference grid, restrict to the set of orthogonal vectors, $\G$, available on the given stencil.  Then the convexified \MA operator~\ref{MAplus} is approximated by
\bq\label{MAmon}\tag*{$(MA)^{M}$}
\mao^{M}[u] \equiv \min\limits_{\{\nu_1\ldots\nu_d\}\in \G}
\prod\limits_{j=1}^{d} \left (\Dt_{\nu_j\nu_j}u\right)^+
\eq

\begin{theorem}[Convergence to Viscosity Solution]\label{thm:viscosity}
Let the PDE~\eqref{MA} have a unique viscosity solution.
Then the solutions of the scheme~\ref{MAmon}, converges to the viscosity solution of~\eqref{MA} as $h,d\theta,\delta\to0$.
\end{theorem}
The proof of convergence follows from verifying consistency and degenerate ellipticity and can be found in~\cite{ABtheory}.

\section{A hybrid discretization}\label{sec:hybrid}\label{sec:disc}
In this section we propose a hybrid discretization of the \MA equation which takes advantage of the best features of each of the previous discretizations. 
 We want to make use of the natural discretization~\ref{MANat} wherever possible in order to take advantage of its simplicity and higher accuracy.  However, we wish to use the monotone discretization~\ref{MAmon} in regions where the solution may be singular in order to properly capture the behaviour of the viscosity solution.  
In this way we hope to achieve the second-order accuracy of the simple discretization in smooth regions and the monotonicity necessary to capture the behaviour of the viscosity solution in non-smooth regions.  

We propose the following hybrid scheme.

Discretize~\eqref{MA} using a weighted average of the two discretizations:
\bq\label{eq:disc}\tag*{$(MA)^H$}
\mao^{H} = w(x)\mao^{N} + (1-w(x))\mao^{M}
\eq
where $w:\Omega\to[0,1]$ is a weight function defined \emph{a priori} from the data as follows.

We first identify $\Omega^s$ which is a neighborhood of the possible singular set of $u$ on $\Omega$, using conditions~\eqref{conditions}. 
\bq\label{SingSet}
\Omega^s = 
\{ x \in \Omega \mid \e < f(x) < 1/\e \}  
\cup  
\{ x \in \partial \Omega \mid \partial \Omega \text{ flat } \text{ or } g(x)  \not \in C^{2,\alpha} \}
\eq
where $\e$ is a small parameter, which we can take equal to $h$, the spatial resolution.

Then define $w(x)$ to be a differentiable function which is zero in an $h$-neighborhood of $\Omega^s$, and which goes to $1$ elsewhere.

\begin{remark} The hybrid scheme will sometimes be less accurate than the standard finite differences when the solution is $C^2$, because it will lose some accuracy at the flat boundary.  While this might seem conservative, there are examples, (see~\cite{BeFrObMA}), where the flat boundary causes blow up in the Hessian, so the monotone scheme is needed.  
\end{remark}

\section{Explicit and semi-implicit solution methods}
Any  discretization of~\eqref{MA} leads to a system of nonlinear equations which must be solved in order to obtain the approximate solution.   

\subsection{Explicit solution methods for monotone schemes}
Using a monotone discretization $F[u]$ of the \MA operator,  the simplest way to solve the \MA equation is by solving the parabolic version of the equation using forward Euler.  That is, we perform the iteration
\[
u^{n+1} = u^n + dt (F[u^n]-f). 
\]
Explicit iterative methods have the advantage that they are simple to implement, but the number of iterations required suffers from the well known CFL condition (which applies in a nonlinear form to monotone discretizations, as explained in~\cite{ObermanSINUM}).  
This approach is slow because for stability it requires a small time step $dt$, which depends on the spatial resolution $h$.  The time step, which can by computed explicitly, is $\bO(h^2)$. This was the approach used in~\cite{ObermanEigenvalues}.
 
\subsection{A semi-implicit solution method}\label{sec:SI}
The next method we discuss is a semi-implicit method, which involves solving the Laplace equation at each iteration.
In~\cite{BeFrObMA} we used the identiy~\eqref{eq:poisson} to build a semi-implicit solution method.  We showed that the method is a contraction, but the strictness of the contraction requires strict positivity of $f$.  In practice, this meant that the iteration was fast for regular solutions, but degenerated to become slower than the explicit method when $f$ was zero in large parts of the domain.

The conditioning of the linearized equation~\eqref{eq:lin}, which affects solution time, 
depends on the strict convexity of the solution, see lemma~\ref{lem:linell}.  The convexity, in turn depends of strict positivity of $f$, see~\autoref{sec:regularity}.
This explains why solution time of the semi-implicit solver depends on regularity. 

Next, we describe a generalization of the semi-implicit method to higher dimensions.  We won't be using the method to solve~\eqref{MA}.  Instead, we will use one iteration to set up the initial value for Newton's method.

Begin with the following identity for the Laplacian in two dimensions,
\bq\label{2d}
\abs{\Delta u} = \sqrt{(\Delta u)^2}
= \sqrt{u_{xx}^2+u_{yy}^2+2u_{xx}u_{yy}}. 
\eq
So if $u$ solves the \MA equation, then
\begin{align*}
\abs{\Delta u} = \sqrt{u_{xx}^2+u_{yy}^2+2u_{xy}^2+2f} = 
\sqrt{ \norm{D^2u}^2 +  2f}
\end{align*}
This leads to a semi-implicit scheme for solving the~\MA equation, used in~\cite{BeFrObMA}. 
\bq\label{eq:poisson}
  \Delta u^{n+1} = \sqrt{2f +  \norm{D^2u^n}^2 } 
\eq
To generalize this to  $\R^d$, we can write the Laplacian in terms of the eigenvalues of the Hessian: $\Delta u = \sum_{i=1}^d \lambda_i[D^2u]$.  
Taking the $d$-th power, and expanding, gives the sum of all possible products of $d$ eigenvalues.
\begin{align*}
(\Delta u)^d
  &= d! \prod\limits_{i=1}^d\lambda_i + P(\lambda_1, \dots, \lambda_d),
\end{align*}
where $P(\lambda)$ is a $d$-homogeneous polynomial, which we won't need explicitly.

The result is the semi-implicit scheme 
\bq\label{SImplicit}
  \Delta u^{n+1} = \sqrt{d! f +   P(\lambda_1[D^2u^n], \dots, \lambda_d[D^2u^n])}. 
\eq
A natural initial value for the iteration is given by the solution of
\bq\label{NMinit}
\Delta u^0 = \sqrt{d! f}.
\eq

\section{Implementation of Newton's method}\label{sec:newton}
To solve the discretized equation
\[
\mao^{H}[u] = f
\]
we use use a damped Newton iteration
\[ u^{n+1} = u^n  - \alpha\delu^n \]
for some $0<\alpha<1$.  The damping parameter $\alpha$ is chosen at each step to ensure that the residual $\|\mao^{H}(u^n)-f\|$ is decreasing.  (In practice we can often take $\alpha = 1$, but damping is sometimes needed.)

The corrector $\delu^n$ solves the linear system
\bq\label{NewtonN2}
\left ( \grad_u \mao^{H} [u^n] \right )\delu^n = \mao^{H}[u^n]-f.
\eq
To set up the equation~\eqref{NewtonN2}, the Jacobian of the scheme is needed.
Since the hybrid discretization is a weighted average of the monotone and standard discretization, and 
 the weight function, $w(x)$, is determined \emph{a priori}, the Jacobian of the hybrid scheme will simply be a weighted average of the corresponding Jacobians.
 
The Jacobian of the \MA operator, discretized using standard finite differences, is given by 
\bq\label{NewtonN}
\grad_u \mao^{N}[u] 
= (\Dt_{xx}u) \Dt_{yy} + (\Dt_{yy}u)\Dt_{xx} - 2 (\Dt_{xy}u)\Dt_{xy},
\eq
which is a discrete version of the linearization of \MA~\eqref{eq:lin}

The Jacobian for the monotone discretization is obtained by using Danskin's Theorem~\cite{Bertsekas} and the product rule.
\[ \grad_u \mao^{M}[u] 
= \sum\limits_{j=1}^d\text{diag}\left( \prod\limits_{k \neq j}\Dt_{\nu_k^*\nu_k^*}u \right)\Dt_{\nu_j^*\nu_j^*}  \]
where the $\nu_j^*$ are the directions active in the minimum in~\ref{MAmon}.

Thus the corrector is obtained by solving the weighted average of the two 
linearizations
\begin{multline}\label{eq:corrector}
 (w(x)\grad_u \mao^{N}[u^n] + (1-w(x)) \grad_u \mao^{M}[u^n])\delu^n \\= w(x)\mao^{N}[u^n] + (1-w(x)) \mao^{M}[u^n]. \end{multline}

In order for the linear equation~\eqref{NewtonN2} to be well-posed, we require
the coefficient matrix to be positive definite. 
As observed in lemma~\autoref{lem:linell}, this condition can fail if the iterate $u^n$ is not strictly convex.

\subsection{Initialization of Newton's method}\label{sec:init}
Newton's method requires a good initialization for the iteration.  Since we need the resulting linear system to be well posed it is essential that the initial iterate: (i) be convex, (ii) respect the boundary conditions, (iii) be close to the solution.   

In order to do this, we first: 
 use one step of the semi-implicit scheme~\eqref{SImplicit}, to obtain a close initial value.  This amounts to solving~\eqref{NMinit} along with consistent Dirichlet boundary conditions~\eqref{Dirichlet}.
Then convexify the result, using the method of~\cite{ObermanCENumerics}.
Since both the steps can be performed on a very coarse grid, and interpolated onto the finer grid, the cost of the initialization is low.

\subsection{Preconditioning}\label{sec:reg}
In degenerate examples, the PDE for $\delu^n$~\eqref{eq:corrector} may be degenerate, which can lead to an ill-conditioned or singular Jacobian.  To get around this problem, we regularize the Jacobian to make sure the linear operator is strictly negative definite; this will not change the fixed points of Newton's method.  We accomplish this by replacing the second directional derivatives $u_{\nu\nu}$ with
\[ \tilde{u}_{\nu\nu} = \max\{u_{\nu\nu},\e\} \]
Here $\e$ is a small parameter.  In the computations of \autoref{sec:2d}, we take $\e = \frac{1}{2dx^2}\ex{8}$.

%
%
%

\section{Computational results in two dimensions}\label{sec:2d} 
In this section, we summarize the results of a number of two-dimensional examples solved using the hybrid scheme described in \autoref{sec:hybrid}.  In particular, we are interested in comparing the computation time for Newton's method with the time required by the methods proposed in~\cite{BeFrObMA}. 
We also visualize the map generated by the gradient of the solution.

These computations are performed on an $N \times N$ grid on the square $[0,1]^2$.  The monotone scheme used a 17 point stencil.  

When needed as part of the initialization, the convex envelope is computed on a coarse grid using the discretization described in~\cite{ObermanCENumerics}.
Since the solution can be computed on a coarse grid, and interpolated, the added computational time is negligible.

\subsection{Failure of Newton's method for natural finite differences}\label{sec:newtonfails}
In this section, we give an example where Newton's method breaks down when standard finite differences are used.   

We chose an example which is only singular at one point, on the boundary.  
Nevertheless, this mild singularity is enough for Newton's method to break down.   

Consider the solution of~\eqref{MA} in $[0,1]^2$, given by 
\[ 
u(\xv) = -\sqrt{2-\norm{\xv}^2}, 
\qquad
f(\xv) = 2 {\left(2-\norm{\xv}^2\right)^{-2} }
\]
The gradient of the solution is unbounded on $ \norm{\xv} = 2$, in particular at the point  $(1,1)$.  The singularity arises from the fact that $f$ is unbounded there.

\begin{figure}[htdp]
\centering
 \subfigure[Solution after two iterations]{\includegraphics[width=.49\textwidth]{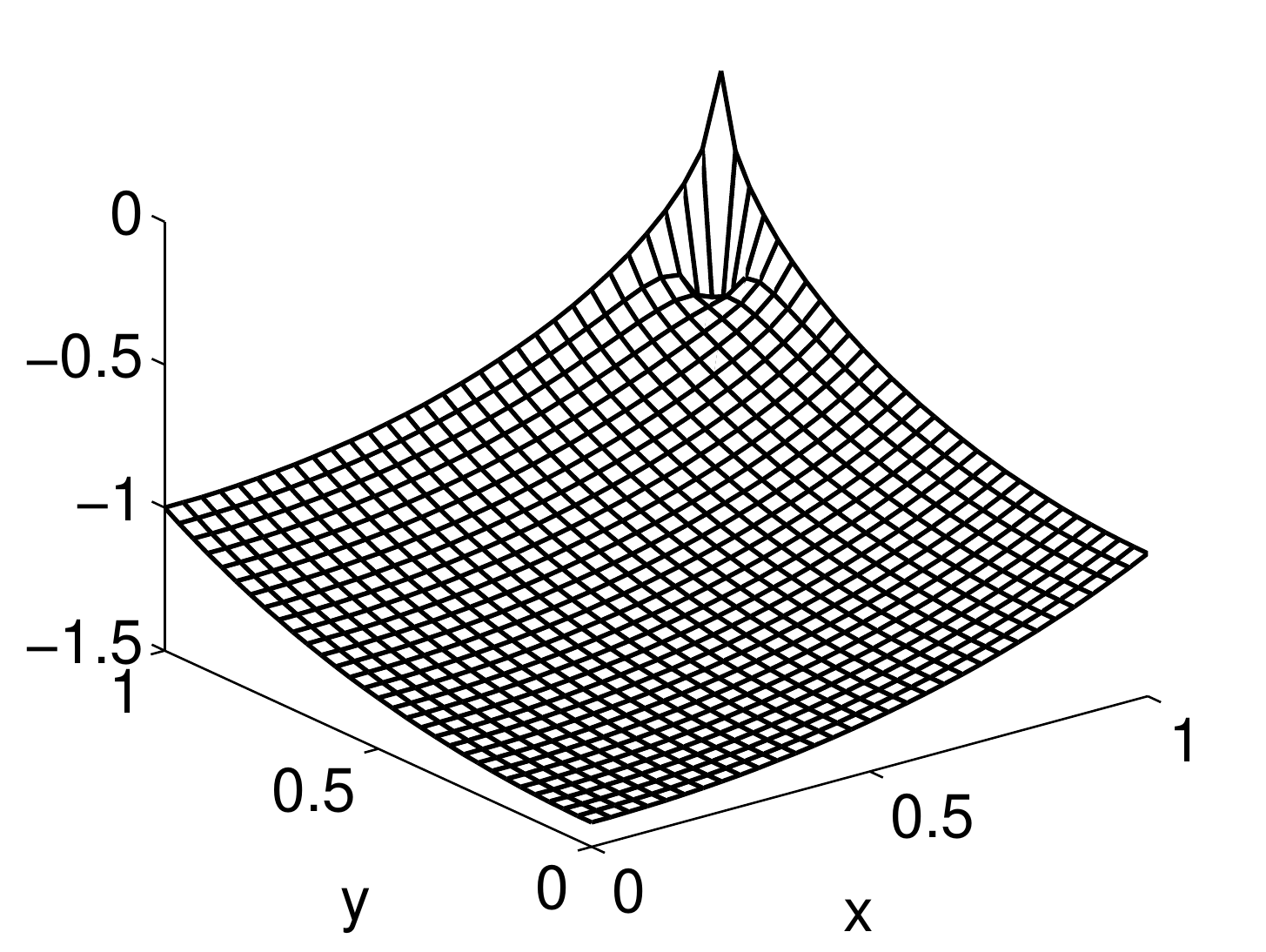}\label{fig:badsol}} 
\subfigure[Gradient map after two iterations]{\includegraphics[width=.49\textwidth]{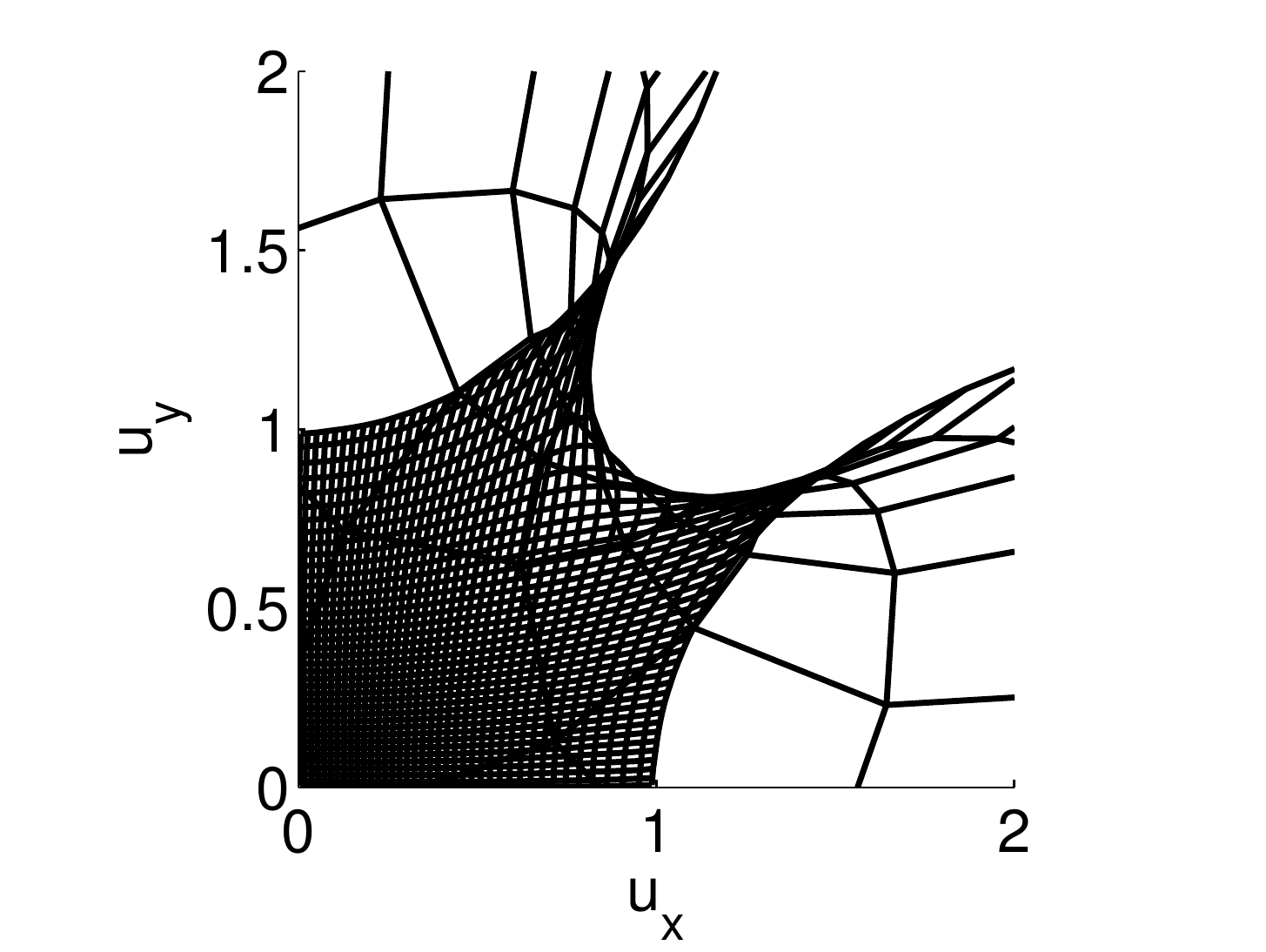}\label{fig:badmesh}}
 \caption{Failure of Newton's method using standard finite differences:
the solution oscillates and becomes non-convex.
}
  \label{fig:fail}
\end{figure}

\begin{figure}
	\centering
	\subfigure[]{\includegraphics[width=.38\textwidth]{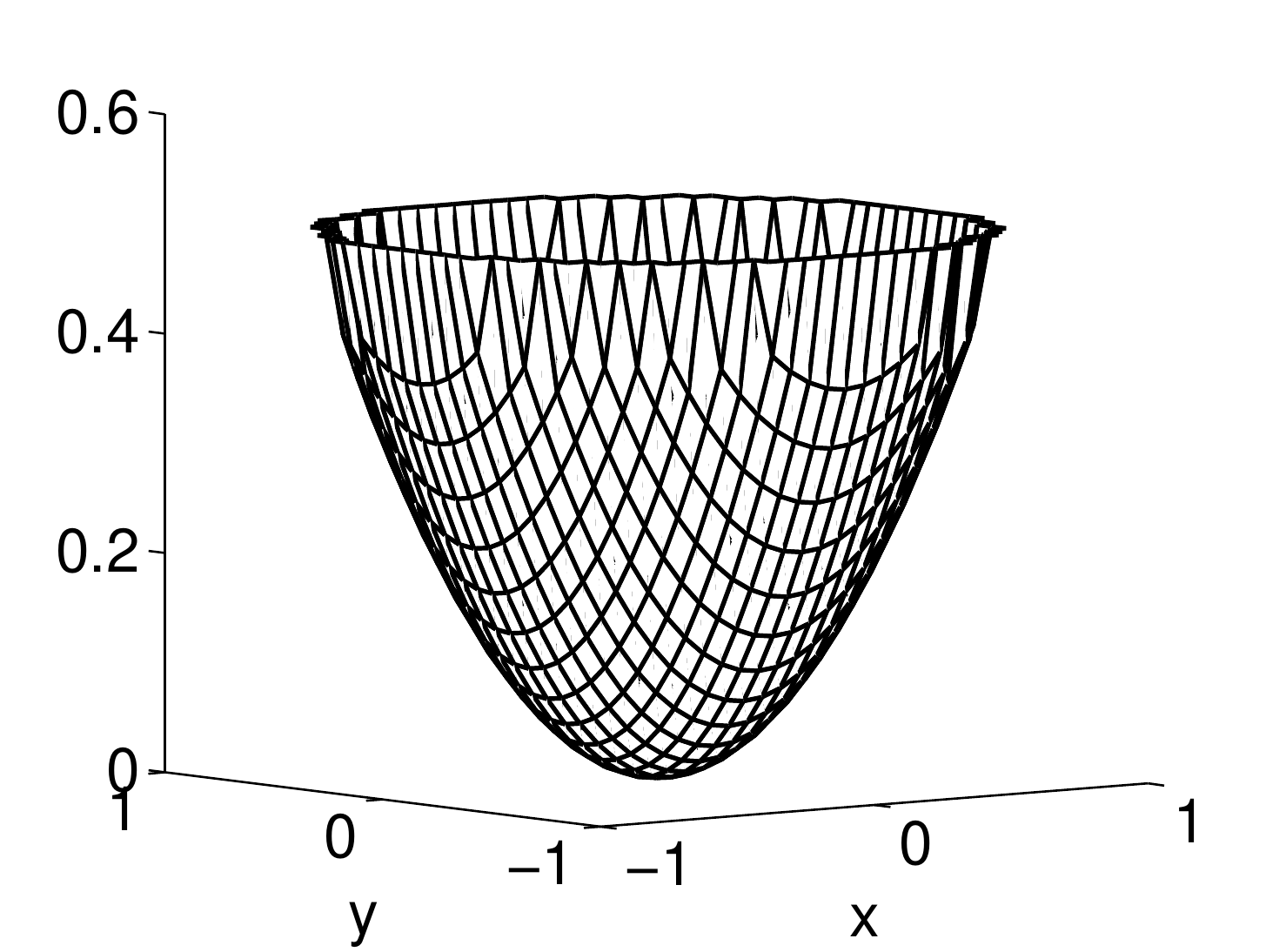}\label{fig:init_surf}}
        \subfigure[]{\includegraphics[width=.38\textwidth]{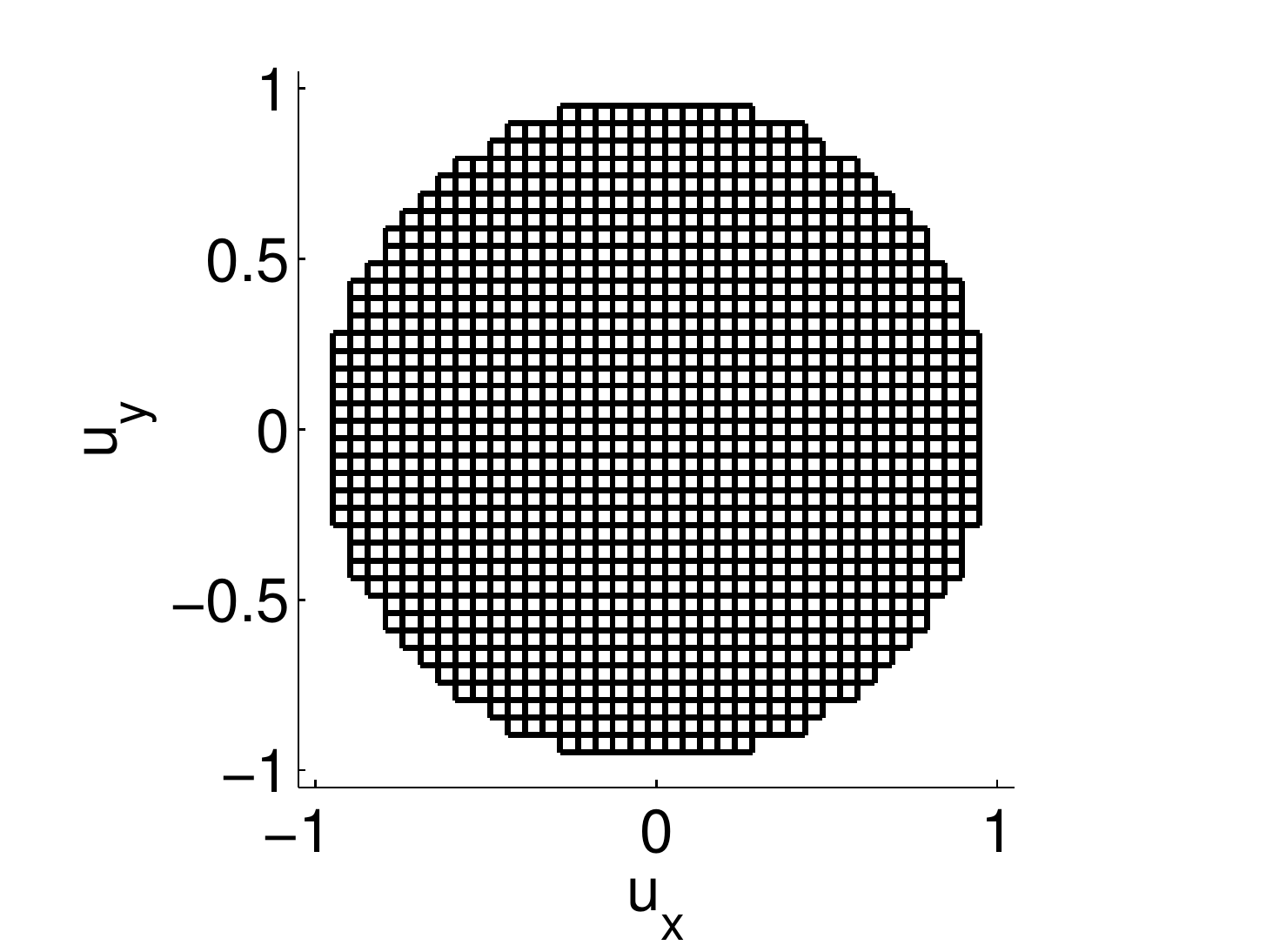}\label{fig:init_mesh}}
        \subfigure[]{\includegraphics[width=.38\textwidth]{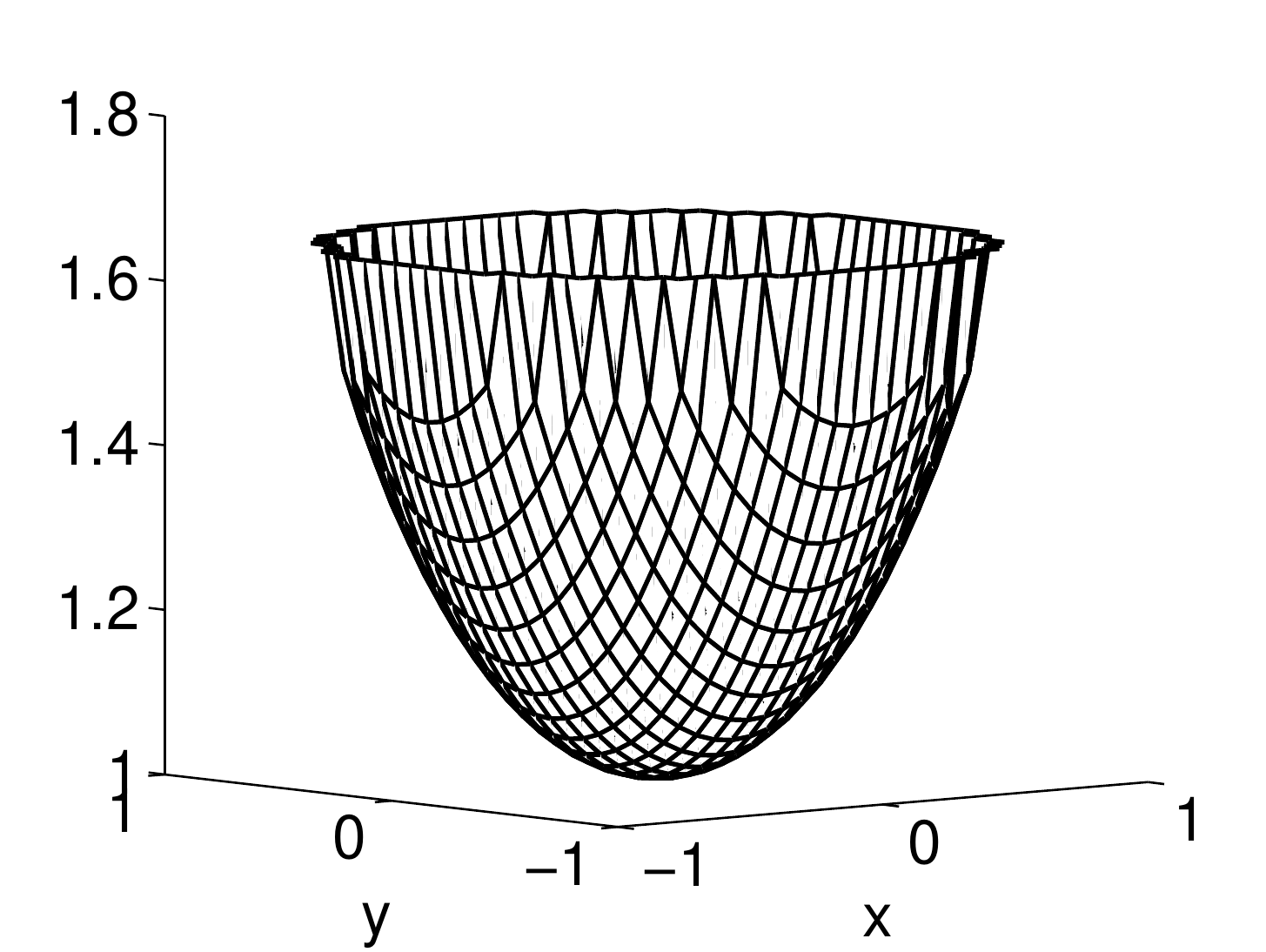}\label{fig:radial_surf}}
        \subfigure[]{\includegraphics[width=.38\textwidth]{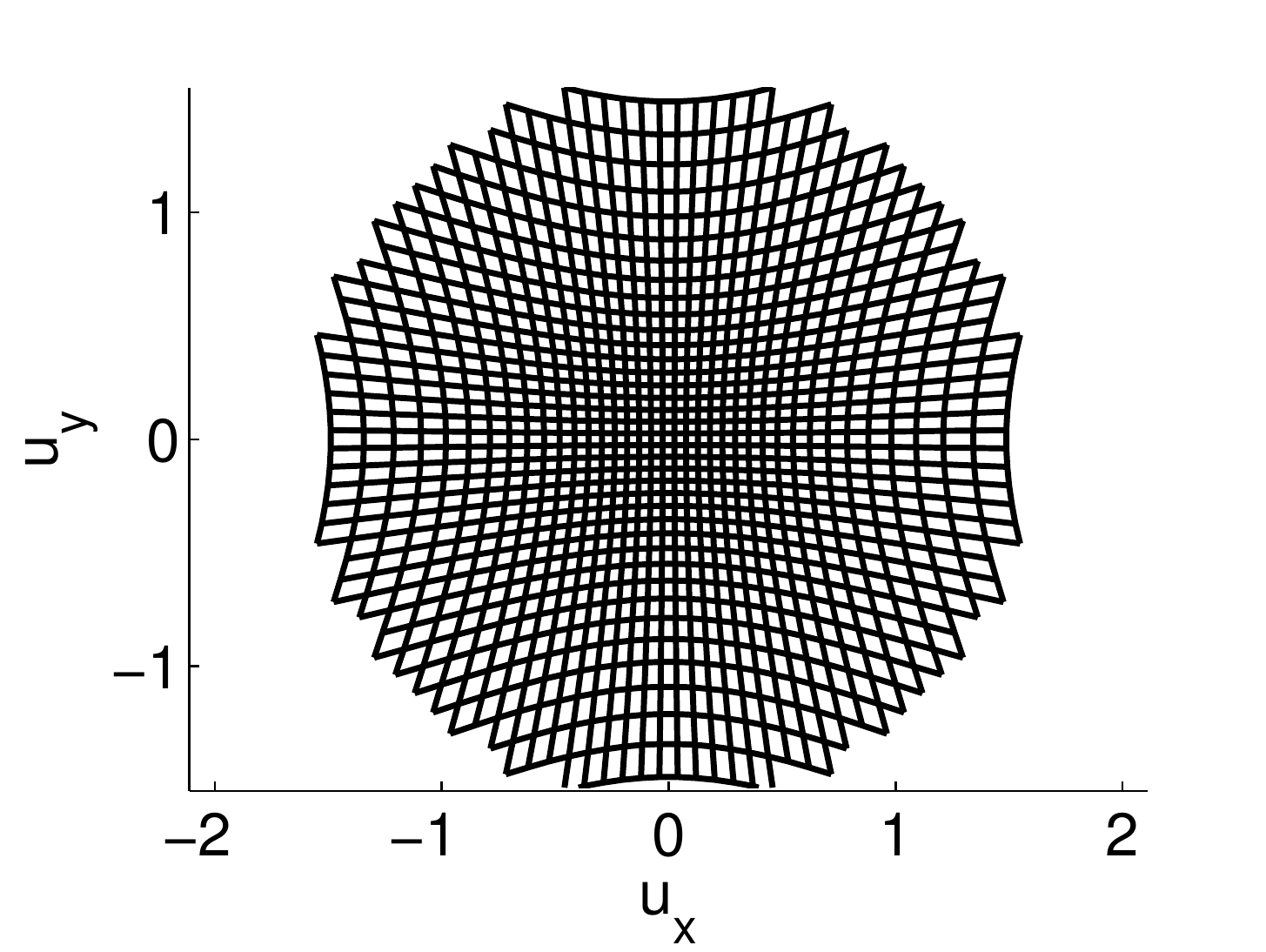}\label{fig:radial_mesh}}	
        \subfigure[]{\includegraphics[width=.38\textwidth]{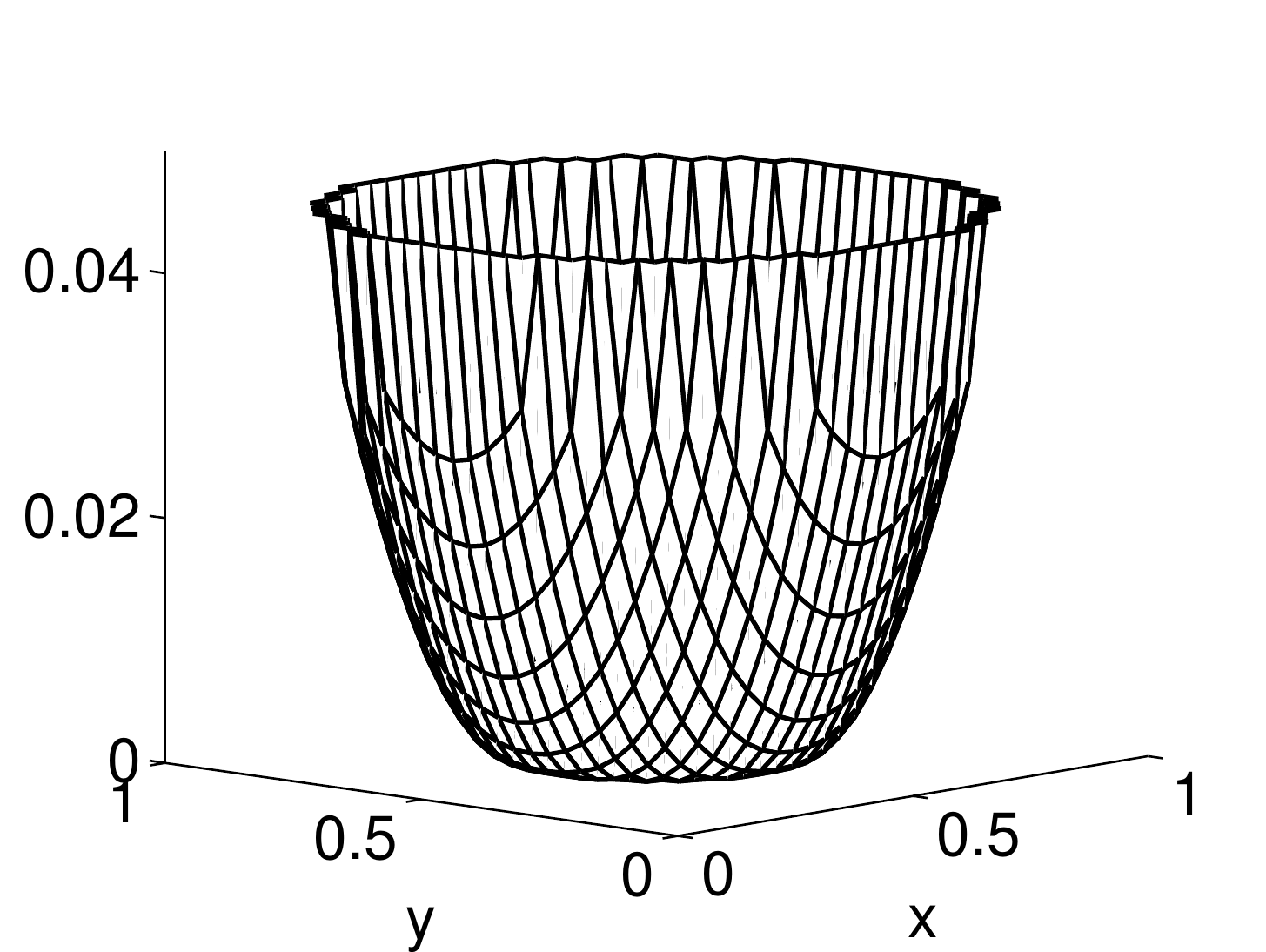}\label{fig:bowl2_surf}}
        \subfigure[]{\includegraphics[width=.38\textwidth]{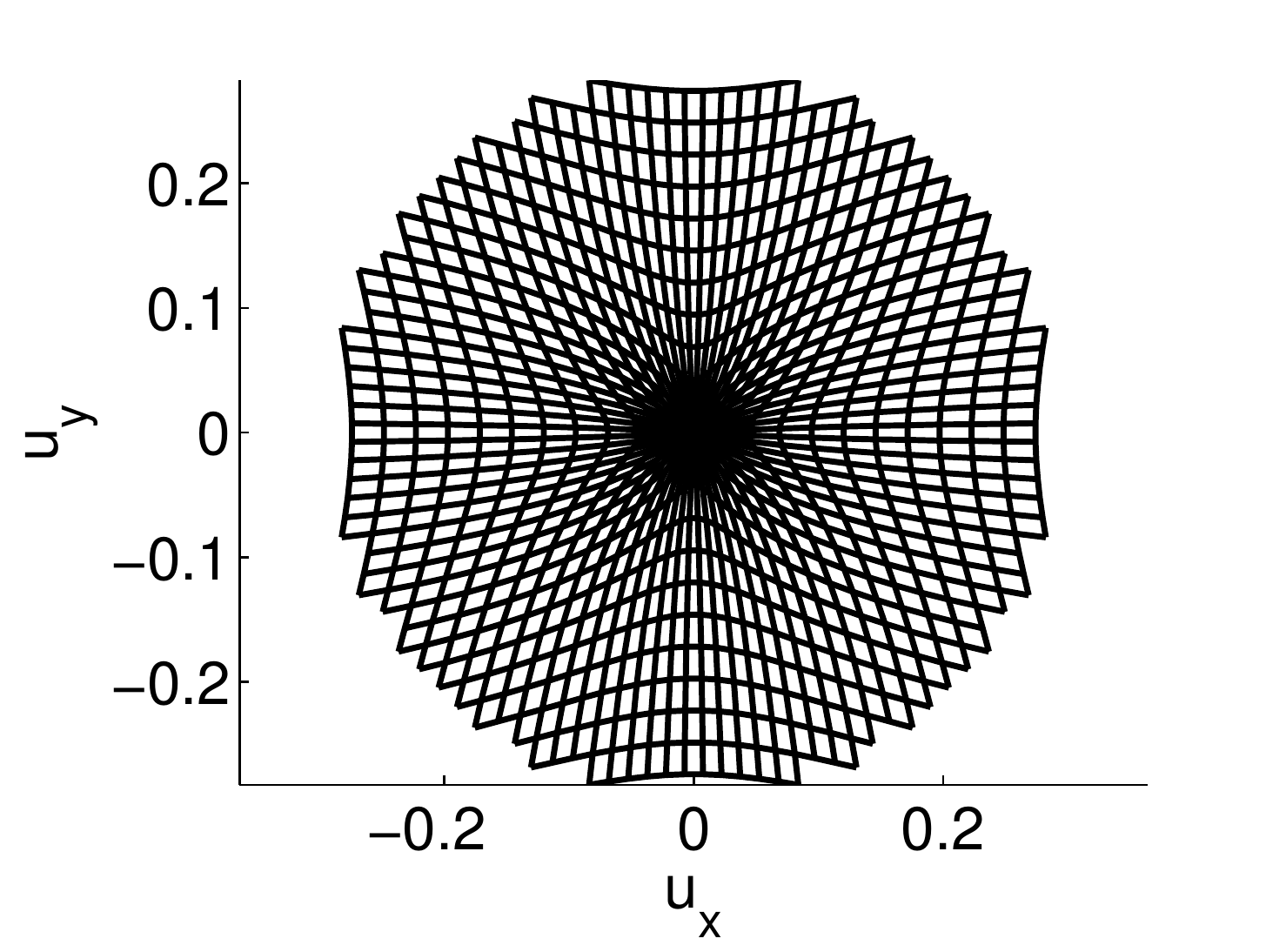}\label{fig:bowl2_mesh}}	
        \subfigure[]{\includegraphics[width=.38\textwidth]{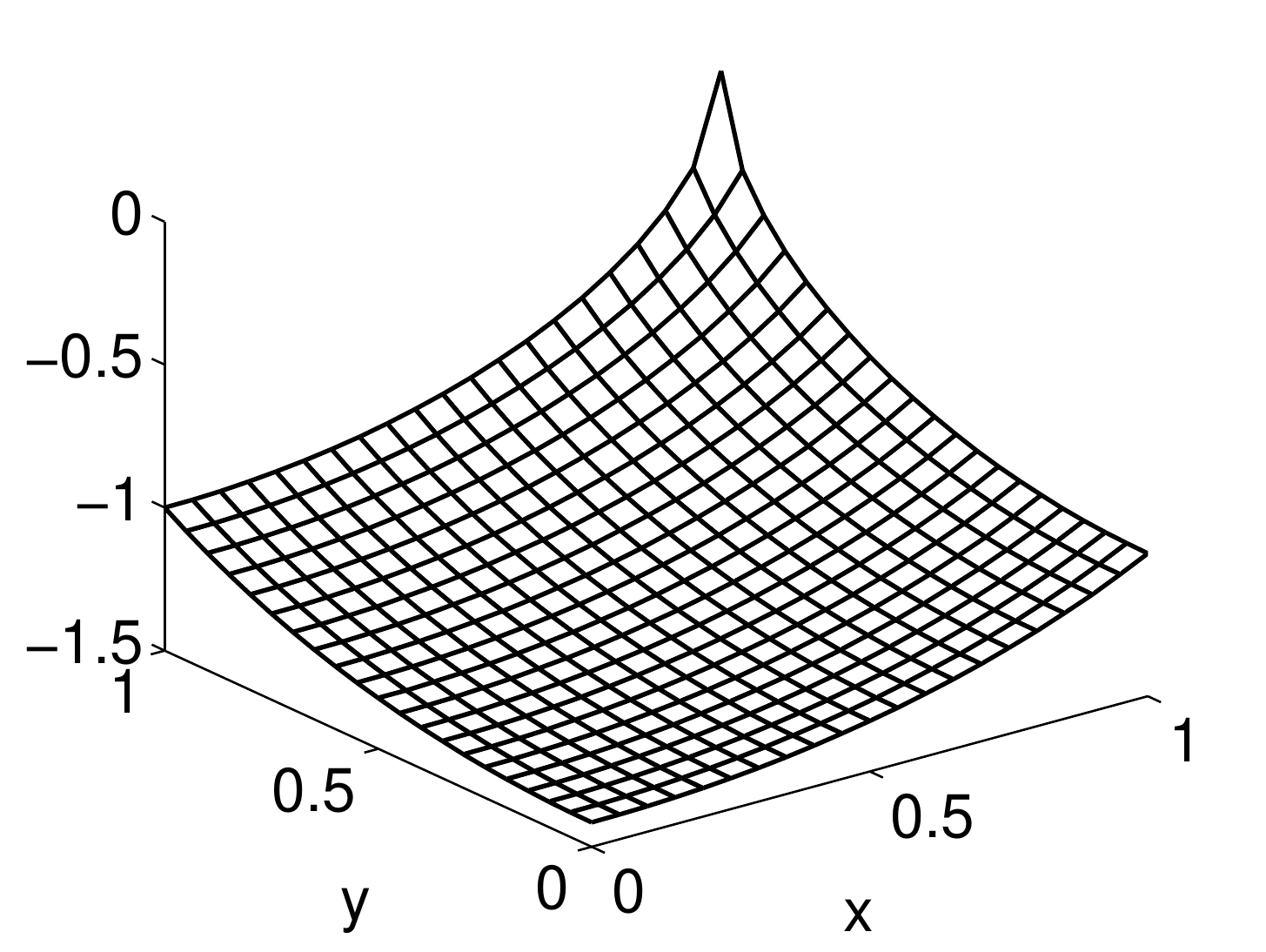}\label{fig:ball_surf}}
        \subfigure[]{\includegraphics[width=.38\textwidth]{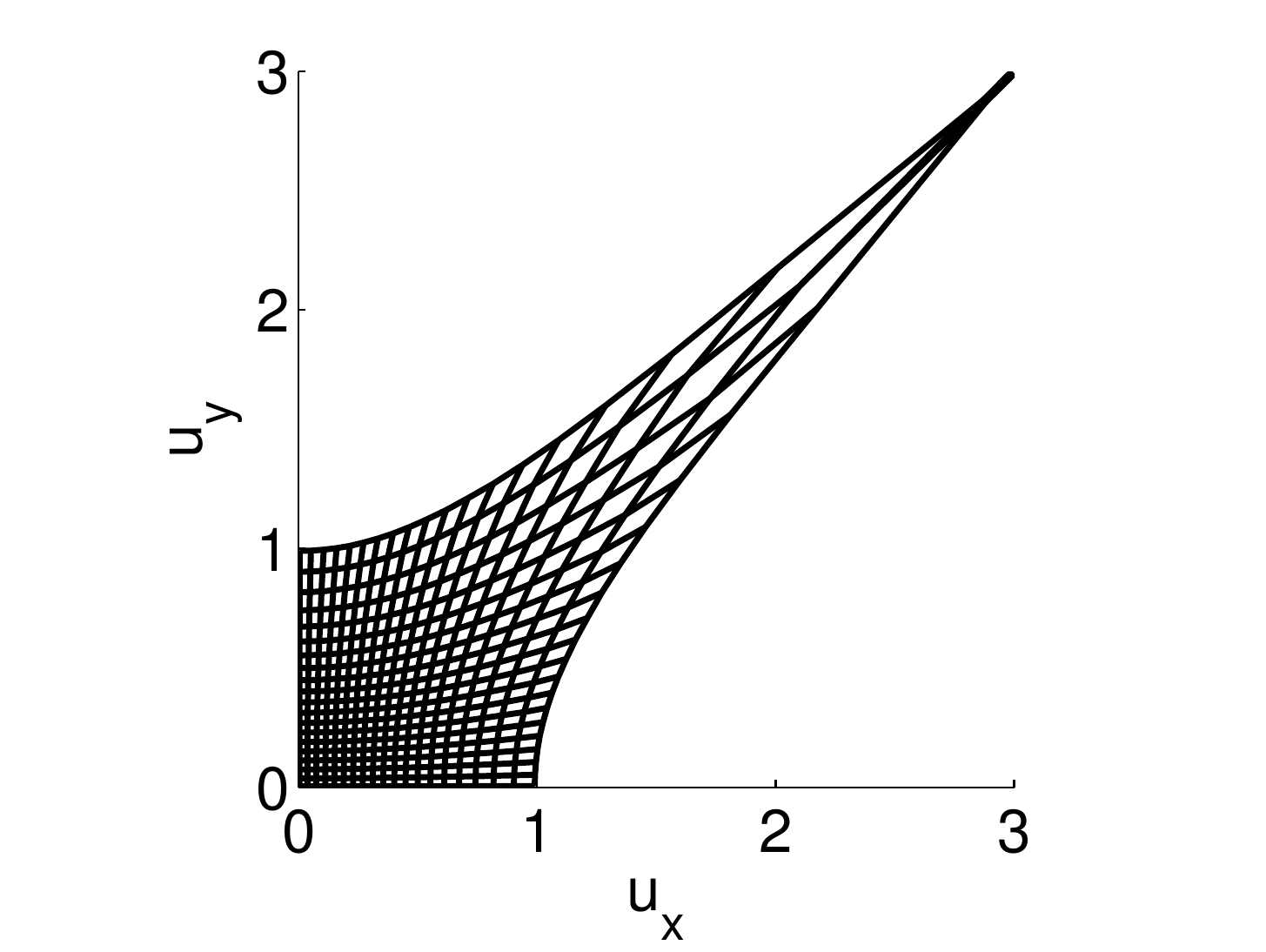}\label{fig:ball_mesh}}
  	\caption{Solutions and mappings for the 
	\subref{fig:init_surf},\subref{fig:init_mesh} identity map, \subref{fig:radial_surf},\subref{fig:radial_mesh} $C^2$ example, 
	\subref{fig:bowl2_surf},\subref{fig:bowl2_mesh} $C^1$ example, and \subref{fig:ball_surf},\subref{fig:ball_mesh} example with blow-up.
	}
  	\label{fig:solutions}
\end{figure}

Due to the singularity, there is an instability in Newton's method if the natural finite difference method is used.
The iteration is initialized with the exact solution.  The result after performing two iterations of Newton's method along with the gradient map, is illustrated in~\autoref{fig:fail}.
The correct computed solution is presented in Figure~\ref{fig:ball_surf}-\ref{fig:ball_mesh}.

\subsection{Four representative examples}
We have tested the hybrid method on a number of examples of varying regularity; the results are summarized in~\autoref{sec:time2d}-\ref{sec:maps}.  To illustrate these results, we present more detailed results for four representative examples.

Write $\xv = (x,y)$, and $\xv_0 = (.5, .5)$ for the center of the domain.

The first example solution, which is smooth and radial, is given by
\bq\label{eq:c2} 
u(\xv) = \exp \left( \frac{ \norm{\xv}^2}{2} \right),
\qquad 
f(\xv) = (1+ \norm{\xv}^2)\exp( \norm{\xv}^2 ).
\eq
The second example, which is $C^1$, is given by
\bq\label{eq:c1} 
u(\xv) = \frac{1}{2}\left( (\norm{\xv-\xv_0} -0.2)^+\right )^2, 
\quad
f(\xv) = 
\left( 
1 - \frac{0.2}{\norm{\xv-\xv_0}}
\right)^+.
\eq

The third example is the one which was used in~\autoref{sec:newtonfails} to demonstrate that Newton's method for standard finite differences is unstable.  
The solution is twice differentiable in the interior of the domain, but has an unbounded gradient near the boundary point $(1,1)$.  The solution is given by
\bq\label{eq:blowup} 
u(\xv) = -\sqrt{2-\norm{\xv}^2},
\qquad 
f(\xv) = 2 {\left(2-\norm{\xv}^2\right)^{-2} }.
\eq

This final is example solution is the cone, which was discussed in~\autoref{sec:alex}.  It is Lipschitz continuous.
\bq\label{eq:cone} 
u(\xv) = \sqrt{\norm{\xv-\xv_0}},
\qquad
f = \mu = \pi \,\delta_{\xv_0}
\eq
It order to approximate the solution on a grid with spatial resolution $h$, using viscosity solutions, we approximate the measure $\mu$ by its average over the ball of radius $h/2$, which gives
\[ 
f^h  = 
\begin{cases}
4/h^2 &  \text{ for } \norm{\xv - \xv_0} \leq h/2,\\
0 & \text{ otherwise.}
\end{cases}
\]

\subsection{Visualization of solutions and gradient maps}\label{sec:maps}
In~\autoref{fig:solutions} the solutions and the gradient maps for the three representative examples are presented. For example~\eqref{eq:cone} the gradient map is too singular to illustrate.
To visualize the maps, the image of a Cartesian mesh under the mapping 
\[ 
\left(\begin{array}{c}x\\y\end{array}\right) \to \left(\begin{array}{c}\Dt_{x}{u}\\ \Dt_{y}{u}\end{array}\right) \] 
is shown, where $(\Dt_x u,\Dt_y u)$ is the numerical gradient of the solution of the \MA equation.  
The image of a circle is plotted for visualization purposes, 
the equation is solved on a square. For reference, the identity mapping is also displayed.  

In each case, the maps agree with the maps obtained using the gradient of the exact solution.

\subsection{Computation time}\label{sec:time2d}
The computation times for the four representative examples is presented in \autoref{table:time2d}.
The computations time are compared to those for the Gauss-Seidel and Poisson iterations described in~\cite{BeFrObMA}.  
The Newton solver is faster in terms of absolute solution time in each case. 
\autoref{table:time} presents of order of magnitude solutions times.
The order of magnitude solution time for Newton's method 
is independent of the regularity of the solutions and faster than both of the other methods.

\begin{table}[htdp]\small
\begin{center}
\begin{tabular}{ccccc}
\multicolumn{5}{c}{$C^2$ Example \eqref{eq:c2}}\\
N & Newton       & \multicolumn{3}{c}{CPU Time (seconds)} \\
  & Iterations   & Newton   &   Poisson  & Gauss-Seidel\\
\hline
31  & 3  & 0.2 & 0.7 & 2.2\\
45  & 4  & 0.2 & 1.1 & 4.1\\
63  & 4  & 0.4 & 1.9 & 15.0\\
89  & 4  & 1.0 & 4.8 & 57.6\\
127 & 5  & 2.9 & 9.6 & 236.7\\
181 & 5  & 9.0 & 23.2 & 1004.0\\
255 & 5  & 30.5 & 52.6 & ---\\
361 & 6  & 131.4 & 162.6 & ---\\
\hline\hline\\
\multicolumn{5}{c}{$C^1$ Example \eqref{eq:c1}}\\
N & Newton       & \multicolumn{3}{c}{CPU Time (seconds)} \\
  & Iterations   & Newton   &   Poisson  & Gauss-Seidel\\
\hline
31  & 4 &  0.4 & 1.1 & 0.8\\
45  & 6  & 0.4 & 6.1 & 2.8\\
63  & 7 & 0.8 & 20.5 & 9.5\\
89  & 9  & 2.0 & 80.0 & 35.9\\
127 & 11  & 5.7 & 256.8 & 145.5\\
181 & 13  & 17.7 & --- & 558.0\\
255 & 16  & 55.3 & --- & ---\\
361 & 20  & 200.0 & --- & ---\\
\hline\hline\\
\multicolumn{5}{c}{Example with blow-up \eqref{eq:blowup}}\\
N & Newton      & \multicolumn{3}{c}{CPU Time (seconds)} \\
  & Iterations   & Newton   &   Poisson  & Gauss-Seidel\\
\hline
31  & 4  &  0.2 & 0.5 & 0.8\\
45  & 4   & 0.4 & 1.4 & 5.3\\
63  & 4  & 0.7 & 2.9 & 19.4\\
89  & 5  & 1.8 & 8.1 & 74.1\\
127 & 7  & 5.1 & 17.7 & 293.3\\
181 & 7  & 12.9 & 51.4 & 1637.1\\
255 & 7  & 36.1 & 128.2 & ---\\
361 & 8  & 152.9 & 374.5 & ---\\
\hline\hline\\
\multicolumn{5}{c}{$C^{0,1}$ (Lipschitz) Example \eqref{eq:cone}}\\
N & Newton       & \multicolumn{3}{c}{CPU Time (seconds)} \\
  & Iterations   & Newton   &   Poisson  & Gauss-Seidel\\
\hline
31  & 9 &  0.5 & 5.3 & 0.8\\
45  & 11  & 0.6 & 27.8 & 5.9\\
63  & 15  & 1.4 & 91.9 & 21.5\\
89  & 22 & 4.3 & 451.0 & 90.5\\
127 & 32  & 14.1 & 1758.2 & 373.9\\
181 & 30  & 34.6 & --- & 1576.1\\
255 & 34 & 101.7 & --- & ---\\
361 & 29  & 280.2 & --- & ---\\
\end{tabular}
\end{center}
\caption{Computation times for the Newton, Poisson, and Gauss-Seidel methods for four representative examples.}
\label{table:time2d}
\end{table}

\begin{table}[htdp]
\begin{tabular}{c|ccc}
 & \multicolumn{3}{c}{Regularity of Solution} \\
Method  
& $C^{2,\alpha}$~\eqref{eq:c2}   
&  $C^{1,\alpha}$~\eqref{eq:c1} and \eqref{eq:blowup}  
& $C^{0,1}$~\eqref{eq:cone}
\vspace{.05cm}\\
\hline
Gauss-Seidel & Moderate  & Moderate &  Moderate \\
             &($\sim\bO(M^{1.8})$) &($\sim\bO(M^{1.9})$) & ($\sim\bO(M^2)$)\\
Poisson      & Fast    & Fast--Slow  & Slow \\
             & ($\sim\bO(M^{1.4}$)  & ($\sim\bO(M^{1.4})$--blow-up) & ($\sim\bO(M^2)$--blow-up)\\
Newton       & Fast     & Fast     & Fast  \\
             & ($\sim\bO(M^{1.3})$) & ($\sim\bO(M^{1.3})$) & ($\sim\bO(M^{1.3})$)\
\vspace{0.2cm}\\
\end{tabular}
\caption{Order of magnitude computation time for the different solvers in terms 
or the regularity of solutions.  Here $M = N^2$ is the total number of grid points.}
\label{table:time}
\end{table}

\subsection{Accuracy}
Numerical errors are presented in~\autoref{table:err2d}.
We compare the accuracy of the hybrid scheme to the standard finite difference discretization, (using the results of~\cite{BeFrObMA}) and to the monotone scheme which was also solved using 
 Newton's method.


We discuss each example in turn.

\subsubsection*{The $C^2$ solution~\eqref{eq:c2}}
The standard finite difference schemes gives $\bO(h^2)$ accuracy (see~\cite{BeFrObMA}).
In this case, the hybrid scheme is slightly \emph{less} accurate, because the monotone scheme is used near the boundary.  
On a strictly convex domain the hybrid scheme would reduce to the standard discretization and achieve the same accuracy.

The effect diminishes as the number of grid points increases so that the number of interior points using the higher order scheme dominates.  
Accuracy approaches $\bO(h^2)$ as the number of grid points increases.  This is a definite improvement over the monotone scheme, which has its accuracy limited by the stencil width.

\subsubsection*{The $C^1$ solution~\eqref{eq:c1}}
The accuracy is $\bO(h)$, which is similar to the standard discretization and better than the limited accuracy permitted by the monotone discretization with a fixed stencil width.  We also look at the error at each point (see~\autoref{fig:error}); it is evident that the singularity around the circle is the factor that most affects the accuracy.  Because of this non-smoothness, there is no reason to expect our scheme to produce the $\bO(h^2)$ accuracy that is possible on $C^2$ solutions.

\subsubsection*{The blow-up solution~\eqref{eq:blowup}}
In this case, the hybrid scheme accuracy is  $\bO(h^{1.5})$.  This is better than the accuracy of both the standard discretization, which was $\bO(h^{0.5})$~\cite{BeFrObMA}, and the monotone scheme, which is limited by the stencil width.

\subsubsection*{The cone solution~\eqref{eq:cone}}
\begin{figure}[htdp]
	\centering
	\subfigure[]{\includegraphics[width=.49\textwidth]{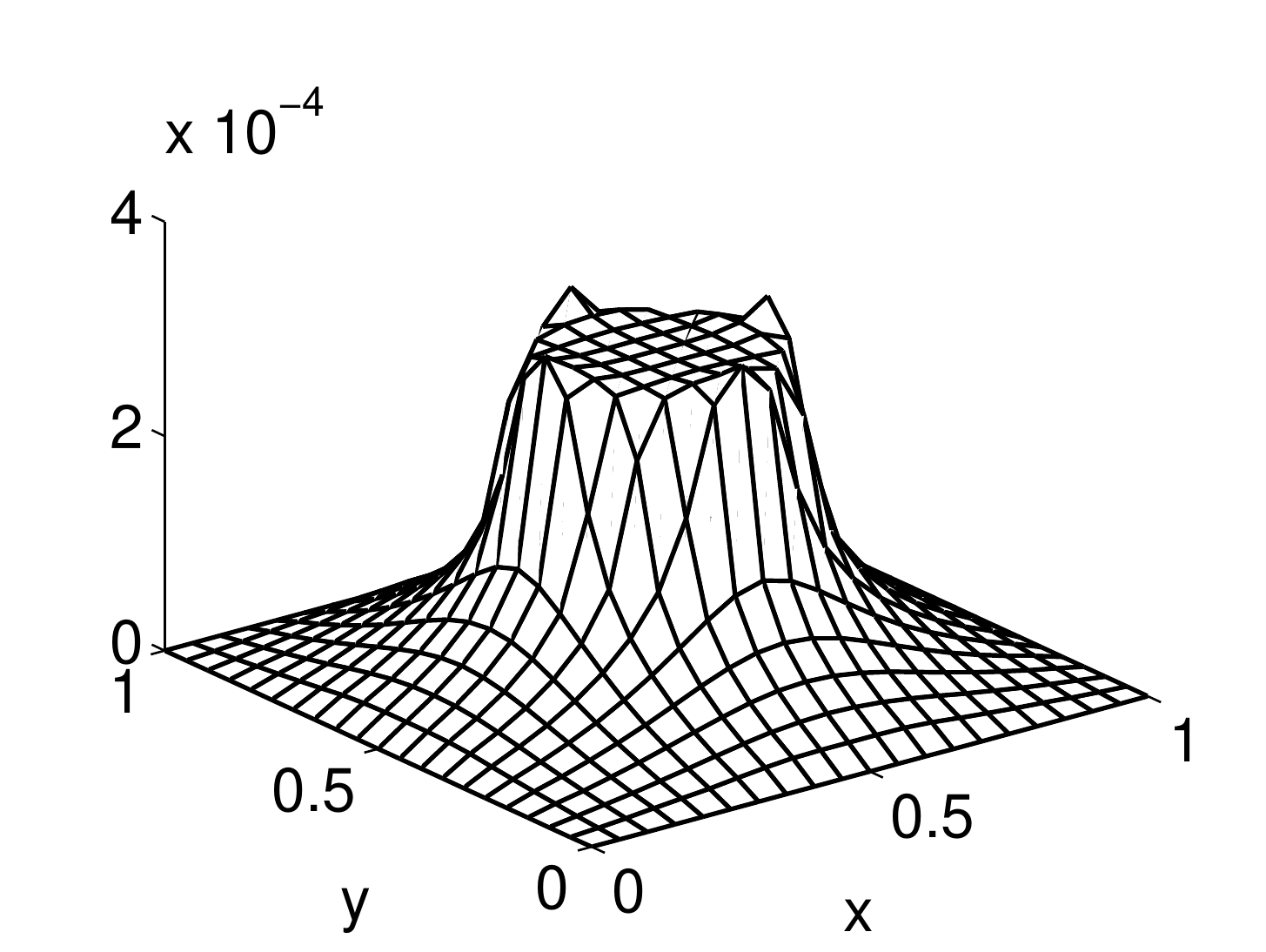}\label{fig:bowlerr}}
        \subfigure[]{\includegraphics[width=.49\textwidth]{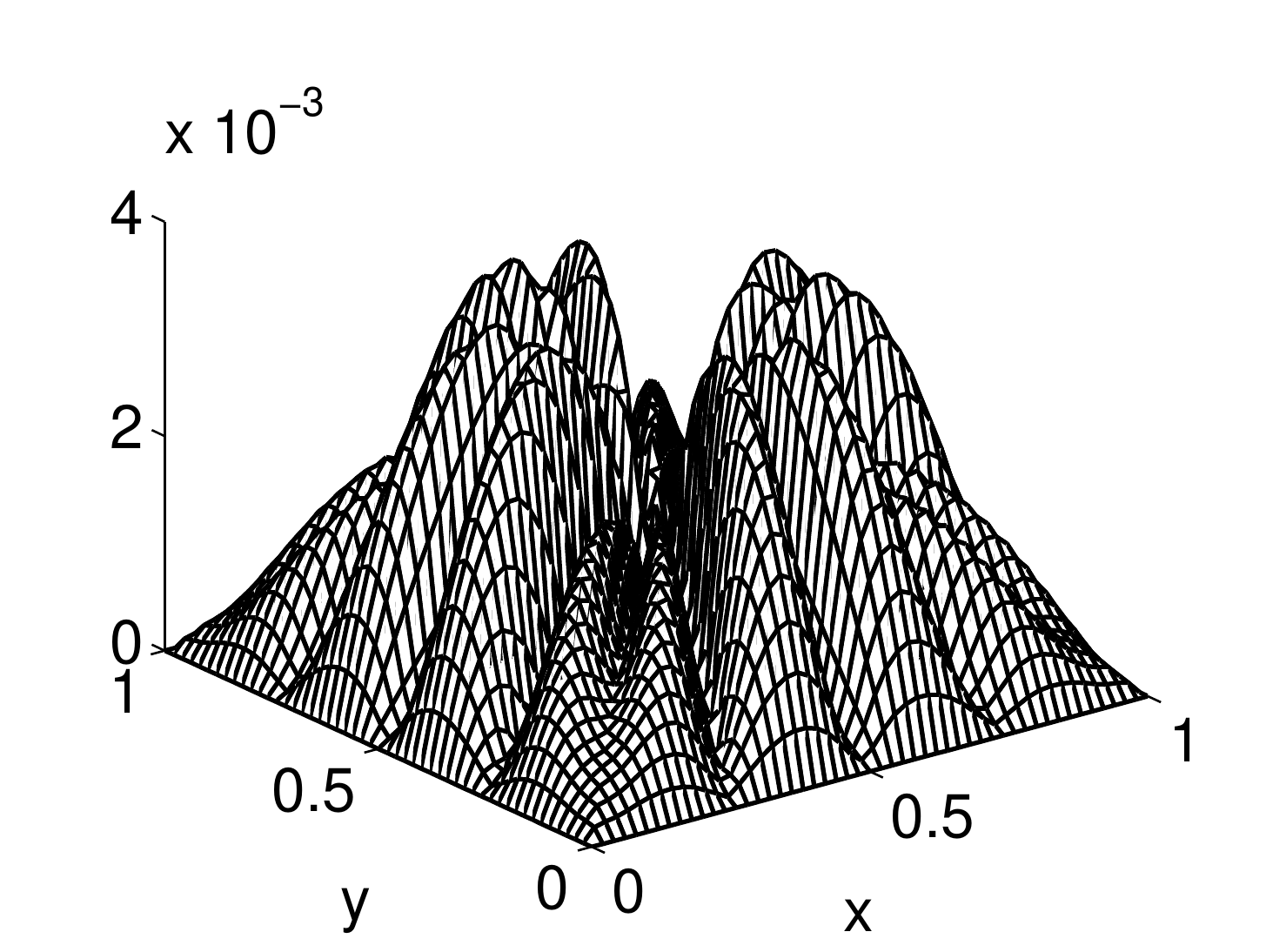}\label{fig:coneerr}}
  	\caption{Surface plots of error using the hybrid scheme for the \subref{fig:bowlerr} $C^1$ example and \subref{fig:coneerr} cone example.}
  	\label{fig:error}
\end{figure}

For this singular example, the hybrid scheme is identical to the monotone scheme (since the right-hand side is either 0 or very large everywhere in the domain).  Consequently, the angular resolution (stencil width) limits the accuracy of solutions.  
We observed that the 17 point stencil reduced the error by an order of magnitude compared to the 9 point stencil.  This dependence on the stencil width is also evident in the surface plot of error (\autoref{fig:error}), which demonstrates that error is largest along directions that are not captured by the stencil.
Since this solution is so singular the reduced accuracy is to be expected.

\begin{table}[htdp]\small
\begin{center}
\begin{tabular}{cccc}
\multicolumn{4}{c}{$C^2$ Example \eqref{eq:c2}}\\
N    & \multicolumn{3}{c}{Maximum Error} \\
    & Standard   &   Monotone  & Hybrid\\
\hline
31  & $7.14\ex{5}$ & $89.09\ex{5}$ & $24.45\ex{5}$\\
45  & $3.39\ex{5}$ & $60.50\ex{5}$ & $15.29\ex{5}$\\
63  & $1.73\ex{5}$ & $50.88\ex{5}$ & $9.06\ex{5}$\\
89  & $0.87\ex{5}$ & $47.51\ex{5}$ & $5.32\ex{5}$\\
127 & $0.43\ex{5}$ & $45.53\ex{5}$ & $3.02\ex{5}$\\
   181 & $0.21\ex{5}$ & $44.65\ex{5}$ & $1.61\ex{5}$\\
255 & $0.11\ex{5}$ & $44.22\ex{5}$ & $0.87\ex{5}$\\
361 & $0.05\ex{5}$ & $44.00\ex{5}$ & $0.46\ex{5}$\\
\hline\hline\\
\multicolumn{4}{c}{$C^1$ Example \eqref{eq:c1}}\\
N    & \multicolumn{3}{c}{Maximum Error} \\
    & Standard   &   Monotone  & Hybrid\\
\hline
31  & $2.6\ex{4}$ & $17.5\ex{4}$ & $12.2\ex{4}$ \\
45  & $1.8\ex{4}$ & $11.6\ex{4}$ & $5.9\ex{4}$ \\
63  & $1.5\ex{4}$ & $9.8\ex{4}$ & $4.2\ex{4}$ \\
89  & $0.9\ex{4}$ & $8.4\ex{4}$ & $2.6\ex{4}$ \\
127 & $0.6\ex{4}$ & $7.9\ex{4}$ & $2.0\ex{4}$ \\
181 & $0.4\ex{4}$ & $7.4\ex{4}$ & $1.2\ex{4}$ \\
255 & --- & $7.2\ex{4}$ & $1.0\ex{4}$ \\
361 & --- & $7.0\ex{4}$ & $0.7\ex{4}$ \\
\hline\hline\\
\multicolumn{4}{c}{Example with blow-up \eqref{eq:blowup}}\\
N    & \multicolumn{3}{c}{Maximum Error} \\
    & Standard   &   Monotone  & Hybrid\\
\hline
31  & $17.15\ex{3}$ & $1.74\ex{3}$ & $1.74\ex{3}$ \\
45  & $14.59\ex{3}$ & $0.98\ex{3}$ & $0.98\ex{3}$ \\
63  & $12.53\ex{3}$ & $0.59\ex{3}$ & $0.59\ex{3}$ \\
89  & $10.67\ex{3}$ & $0.37\ex{3}$ & $0.35\ex{3}$ \\
127 & $9.00\ex{3}$ & $0.35\ex{3}$ & $0.20\ex{3}$ \\
181 & $7.59\ex{3}$ & $0.34\ex{3}$ & $0.12\ex{3}$ \\
255 & $6.42\ex{3}$ & $0.33\ex{3}$ & $0.07\ex{3}$ \\
361 & $5.41\ex{3}$ & $0.33\ex{3}$ & $0.04\ex{3}$ \\
\hline\hline\\
\multicolumn{4}{c}{$C^{0,1}$ (Lipschitz) Example \eqref{eq:cone}}\\
N    & \multicolumn{3}{c}{Maximum Error} \\
    & Standard   &   Monotone  & Hybrid\\
\hline
31  & $10\ex{3}$ & $3\ex{3}$ & $3\ex{3}$ \\
45  & $8\ex{3}$ & $3\ex{3}$ & $3\ex{3}$ \\
63  & $6\ex{3}$ & $3\ex{3}$ & $3\ex{3}$ \\
89  & $4\ex{3}$ & $4\ex{3}$ & $4\ex{3}$ \\
127 & $3\ex{3}$ & $4\ex{3}$ & $4\ex{3}$ \\
181 & $2\ex{3}$ & $4\ex{3}$ & $4\ex{3}$ \\
255 & --- & $4\ex{3}$ & $4\ex{3}$ \\
361 & --- & $4\ex{3}$ & $4\ex{3}$ 
\end{tabular}
\end{center}
\caption{Accuracy for the standard, monotone, and hybrid discretizations for four representative examples.}
\label{table:err2d}
\end{table}

\section{Computational results in three dimensions.}\label{sec:3d}
In this section, we demonstrate the speed and accuracy of the hybrid Newton's method for three dimensional problems.
These computations are performed on an $N \times N \times N$ grid on the square $[0,1]^3$.  The monotone scheme used a 19 point stencil.

The size of the computation was restricted by the available memory, 
 not by solution time (the computations were performed on a recent model laptop).

The solution methods of~\cite{BeFrObMA} were restricted to the two-dimensional \MA equation, so we are no longer able to compare solution times to Newton's method for these examples. 

As before, we provide specific results for three representative examples of varying regularity.  In this section we use the notation
\[ \xv = (x,y,z) \]
and let  $\xo = (.5, .5, .5)$ be the centre of the domain.

The first example is the $C^2$ solution given by
\bq\label{eq:c23d} 
u(\xv) = \exp{\left(\frac{\abs{\xv}^2}{2}\right)}, 
\qquad 
f(\xv) = (1+\abs{\xv}^2)\exp{(\frac{3}{2}\abs{\xv}^2)}. 
\eq
The second example is the $C^1$ solution given by
\bq\label{eq:c13d} 
u(\xv) = \frac{1}{2}\left( (\norm{\xv-\xv_0} -0.2)^+\right )^2, 
\eq
\[ f(\xv) = \begin{cases}
1 - \frac{0.4}{\abs{\xv-\xo}}+\frac{0.04}{\abs{\xv-\xo}^2}, 
 &\abs{\xv-\xo}>0.2\\
0 & \text{otherwise.}
\end{cases}\]
The third example is the surface of a ball, which is differentiable in the interior of the domain, but has an unbounded gradient at the boundary.
\bq\label{eq:blowup3d}u(\xv) = -\sqrt{3-\abs{\xv}^2},\qquad f(\xv) = 3(3-\abs{\xv}^2)^{-5/2}.\eq

As indicated by the results in \autoref{table:3d}, the hybrid Newton's method continues to perform well in three dimensions.
(The fact that the solver required only one iteration for Example~\eqref{eq:c13d} was simply an artifact---for larger problems sizes more iterations were required.

\begin{table}[htdp]
\begin{center}
\begin{tabular}{cccc}
\multicolumn{4}{c}{$C^2$ Example \eqref{eq:c23d}}\\
N & Max Error & Iterations  & CPU Time (s) 
\vspace{.05cm}\\
\hline
7  & 0.0151 & 2 & 0.04 \\
11 & 0.0140 & 3 & 0.10 \\
15 & 0.0129 & 5 & 0.71 \\
21 & 0.0121 & 6 & 6.72 \\
31 & 0.0111 & 5 & 86.63\\
\hline\hline\\
\multicolumn{4}{c}{$C^1$ Example \eqref{eq:c13d}}\\
N & Max Error & Iterations  & CPU Time (s) 
\vspace{.05cm}\\
\hline
7 & 0.0034 & 1 & 0.02\\
11 & 0.0022 & 1 & 0.09\\
15 & 0.0016 & 1 & 0.22\\
21 & 0.0009 & 1 & 1.03\\
31 & 0.0005 & 1 & 17.12\\
\hline\hline\\
\multicolumn{4}{c}{Example with Blow-up \eqref{eq:blowup3d}}\\
N & Max Error & Iterations  & CPU Time (s) 
\vspace{.05cm}\\
\hline
7  &  $9.6\ex{3}$ & 1 & 0.03\\
11 & $5.2\ex{3}$ & 3 & 0.11\\
15 & $4.6\ex{3}$ & 3 & 0.48\\
21 & $4.0\ex{3}$ & 6 & 7.42\\
31 & $2.9\ex{3}$ & 8 & 138.74
\end{tabular}
\end{center}
\caption{Maximum error and computation time for the hybrid Newton's method on three representative examples.}
\label{table:3d}
\end{table}

\section{Conclusions}
The purpose of this work was to build a fast, accurate finite difference solver for the elliptic \MA equation.

A hybrid finite difference discretization was used 
which selects between an accurate standard finite difference discretization and a stable (provably convergent) monotone discretization.  
The choice of discretization was based on known regularity results which depended on the boundary data, $g$, the right hand side function $f$, and strict convexity of the domain.  Wherever the requirements on the data are not met, the hybrid discretization chooses the monotone discretization.

The discretized equations were solved by Newton's method, which is fast, $\bO(M^{1.3})$, where $M$ is the number of data points, independent of the regularity of the solution.
The implementation of Newton's method was significantly (orders of magnitude) faster than the two other methods used for comparison.  
The hybrid discretization was shown to be necessary for stability of Newton's method: an example with a mildly singular solution showed that the standard discretization leads to instabilities.

The hybrid discretization was introduced to improve the accuracy of the monotone discretization on regular solutions.   This expected  improvement was achieved.   On regular solutions the hybrid solver was (asymptotically) as accurate as the standard finite difference discretization.    For one moderately singular example the hybrid solver was more accurate than standard finite differences by $\bO(h)$.

The discretization and solution method used was not restricted to two dimensions.  This allowed for the solution of three dimensional problems on moderate  sized grids, with the problem size limited by computer memory, not solution time.

In summary, the solver presented used a novel discretization in general dimensions, accompanied by a fast solution method.  The resulting solver is a significant improvement over existing methods for the solution of possibly singular solutions of the elliptic \MA equation, in terms of solution time, stability, and accuracy.  


\bibliographystyle{model1-num-names}

\bibliography{MongeAmpere}

\end{document}